\documentclass[final,3p,times]{elsarticle}
\usepackage{amssymb}
\usepackage{amsmath}
\usepackage{algorithm}
\usepackage{algpseudocode}
\usepackage{booktabs}
\usepackage{graphicx} % Required for inserting images
\usepackage{indentfirst}
\usepackage{mathtools}
\usepackage{multirow}
\usepackage{subcaption}

\newcommand{\onehalf}{\frac{1}{2}}
\DeclarePairedDelimiter\ceil{\lceil}{\rceil}

\journal{Journal of Computational Physics}

\begin{document}

\begin{frontmatter}


\title{Comparison between \textit{a priori} and \textit{a posteriori} slope limiters for high-order finite volume schemes}

\author[label1]{Jonathan Palafoutas}
\author[label2]{David A. Velasco Romero}
\author[label1,label2]{Romain Teyssier}
\affiliation[label1]{organization={Program in Applied \& Computational Mathematics, Princeton University},
         addressline={Washington Road},
         city={Princeton},
         postcode={08544},
         state={NJ},
         country={USA}}

\affiliation[label2]{organization={Department of Astrophysical Sciences, Princeton University},
         addressline={171 Broadmead Street},
         city={Princeton},
         postcode={08540},
         state={NJ},
         country={USA}}

%% Abstract
\begin{abstract}

High-order finite volume and finite element methods offer impressive accuracy and cost efficiency when solving hyperbolic conservation laws with smooth solutions. However, if the solution contains discontinuities, these high-order methods can introduce unphysical oscillations and severe overshoots/undershoots. Slope limiters are an effective remedy, combating these oscillations by preserving monotonicity. Some limiters can even maintain a strict maximum principle in the numerical solution. They can be classified into one of two categories: \textit{a priori} and \textit{a posteriori} limiters. The former revises the high-order solution based only on data at the current time $t^n$, while the latter involves computing a candidate solution at $t^{n+1}$ and iteratively recomputing it until some conditions are satisfied. These two limiting paradigms are available for both finite volume and finite element methods.

In this work, we develop a methodology to compare \textit{a priori} and \textit{a posteriori} limiters for finite volume solvers at arbitrarily high order. We select the maximum principle preserving scheme presented in \cite{zhang2011maximum, zhang2010maximum} as our \textit{a priori} limited scheme. For \textit{a posteriori} limiting, we adopt the methodology presented in \cite{clain2011high} and search for so-called \textit{troubled cells} in the candidate solution. We revise them with a robust MUSCL fallback scheme. The linear advection equation is solved in both one and two dimensions and we compare variations of these limited schemes based on their ability to maintain a maximum principle, solution quality over long time integration and computational cost.

This analysis reveals a fundamental tradeoff between these three aspects. The high-order \textit{a posteriori} limited solutions boast great quality at long time-scales, taking full advantage of the sharp gradients of the high-order finite volume method. However, they introduce consistent maximum principle violations. On the other hand, the high-order \textit{a priori} limited solutions can preserve a strict maximum principle. Interestingly, this is still true when the classic fourth-order Runge-Kutta method is used, despite it not being classified as strong-stability-preserving. However, we find a serious drawback; with a spatial polynomial reconstruction of degree five or higher, the \textit{a priori} limited solution becomes dominated by numerical artifacts and diffusion, exhibiting much worse solution quality than their \textit{a posteriori} counterparts. Convex blending of revised fluxes reduces the magnitudes of the violations produced by \textit{a posteriori} limited schemes, but at the expense of more visible numerical artifacts.

Moreover, we compare two methods for computing the flux integrals along two-dimensional cell faces, revealing that one option is more cost-effective but leads to particularly large violations when used with the \textit{a priori} limited scheme. Consequently, the \textit{a priori} limited scheme is forced to use the more costly flux computation, making it significantly more expensive at higher-order than the \textit{a posteriori} limited scheme on CPU architecture. This cost difference can be almost entirely mitigated using a GPU implementation of the same schemes, highlighting that GPUs are well-suited for high-order finite volume stencil operations.

\end{abstract}

\begin{graphicalabstract}
\end{graphicalabstract}

\begin{highlights}
\item Developed experimental methodology with which to compare \textit{a priori} and \textit{a posteriori} slope limiters for arbitrarily-high-order finite volume schemes.
\item Conducted numerical tests for the linear advection equation in one and two dimensions.
\item Demonstrated a tradeoff between maximum principle violations, long time-scale numerical solution quality, and computational cost.
\end{highlights}

%% Keywords

\end{frontmatter}



\section{Introduction}

There have been significant advancements in the development of numerical schemes for hyperbolic conservation laws over the past few decades. With increasing demands for precision and reliability in simulations across fields like fluid dynamics, astrophysics, and climate modeling, there is a growing trend toward conducting simulations at higher resolutions.

More recently, there has been a shift in focus towards high-order numerical schemes, typically referring to those of order three or higher. These schemes offer several advantages over their lower-order counterparts, including improved resolution of fine details, reduced numerical diffusion over extended time scales, and the capability to capture sharper gradients and discontinuities. Moreover, in problems with smooth solutions, higher-order schemes exhibit exponentially lower errors as the order increases, meaning they can achieve the same accuracy as a given low-order scheme with a significantly lower resolution. 

However, they also come with a significant drawback when dealing with solutions that contain discontinuities:  Using higher-degree interpolation polynomials can create unphysical oscillations and severely overshoot/undershoot the true solution bounds. In problems involving the conservation of a scalar $u(t, \mathbf{x})$, it may be essential for the numerical solution to satisfy the maximum principle, which defines a range of values represented by $M=\max_\mathbf{x}u_0(\mathbf{x})$ and $m=\min_\mathbf{x}u_0(\mathbf{x})$, where $u_0(\mathbf{x})=u(0,\mathbf{x})$ is the initial condition. These spurious oscillations violate the maximum principle condition $u(t, \mathbf{x}) \in [m, M]$,  making it challenging to enforce physical properties such as the positivity of density or pressure. In fact, Godunov's theorem asserts that any linear scheme higher than first-order will not strictly preserve the maximum principle of a scalar conservation law \cite{godunov1959finite}. Maximum principle violations are catastrophic for solvers of many nonlinear conservation laws, rendering linear high-order schemes useless for such problems.

As far as linear schemes for conservation laws go, there are three popular spatial discretizations: finite difference (FD), finite volume (FV), and finite element (FE) methods. FD/FV methods can be made higher-order simply by increasing the size of the stencil used to interpolate spatial derivatives from nodes/cells and their neighbors. FV methods are particularly attractive for conservation laws because they can easily be defined to conserve solution quantities to numerical precision. Additionally, as demonstrated in the text, they remain stable with explicit Runge-Kutta methods even for large CFl factors.

In the realm of scientific computing, finite element (FE) methods, especially the discontinuous Galerkin (DG) method, have gained increasing popularity.  In DG, the numerical solution is represented using a polynomial basis within a domain region known as an element. By using higher-degree polynomials across fewer elements, the number of discontinuities between solutions is relatively low compared to an equivalent FV scheme. This aspect is crucial, particularly when addressing these discontinuities involves costly Riemann solvers.

Many techniques have been developed to help numerical solutions preserve local monotonicity and avoid these violations. For instance, the artificial viscosity method adds a physical diffusion term to the underlying equation in order to smooth numerical solutions near discontinuities or elsewhere \cite{vonneumann1950method, donnert2019weno, colella1984piecewise, premasuthan2014computation, lu2019nonlinear}.

Another approach is to utilize a slope limiter, an operator designed to reduce the slopes of the piecewise reconstruction of a solution, thereby enforcing monotonicity. The key point is that this operator is nonlinear, allowing schemes to break free from the constraints of Godunov's theorem. There are many examples of limited schemes in the literature. The MUSCL scheme introduced by Van Leer \cite{van1979towards} is a second-order scheme that strictly preserves a given maximum principle. Later came the popular TVD limiter \cite{sweby1984high} and the PPM scheme \cite{sweby1984high}. The ENO \cite{harten1997uniformly} scheme and WENO scheme \cite{hu1999weighted} were introduced as a means to preserve the high-order accuracy of FV schemes  while reducing oscillations near discontinuities.

A strictly maximum principle-preserving (MPP) scheme working at arbitrarily high order was introduced by Zhang \& Shu \cite{zhang2011maximum, zhang2010maximum}. Their method calculates the numerical solution as a convex combination of high-order and first-order updates. It is guaranteed to satisfy the maximum principle if a sufficiently small CFL factor is used along with first-order forward Euler integration or another Runge-Kutta method equivalent to a convex combination of Euler steps. These methods are known as Strong-Stability-Preserving (SSP) integrators. While the popular fourth-order explicit Runge-Kutta method is not classified as Strong Stability Preserving (SSP), it exhibits quasi-SSP behavior for certain problems \cite{sanz2010positivity}. Therefore, we will evaluate its compatibility with MPP schemes.

it behaves as a \textit{quasi}-SSP integrator, enabling schemes to maintain a strict maximum principle.

While slope limiters are popular tools when solving discontinuous problems with FV schemes, and can even guarantee the avoidance of maximum principe violations, they have drawbacks for FE. These popular limiters often necessitate reducing the degree of interpolating polynomials to zero to prevent overshoots and undershoots. However, in FE methods where polynomial elements cover large spatial regions, reducing an element to first-order significantly compromises the accuracy of the overall solution.

An important advancement in addressing this issue of slope limiting for FE was made by considering \textit{a posteriori} limiters. Unlike traditional \textit{a priori} limiters, which use data at $t^n$ to compute a limited update at $t^{n+1}$ in explicit schemes, \textit{a posteriori} limiters compute a candidate solution at $t^{n+1}$ and revise it until it satisfies the problem requirements. This concept, introduced by Krivodonova \cite{krivodonova2007limiters} and Clain \textit{et al.} \citep{clain2011high} as MOOD, involves reducing an element solution by one degree at a time, limiting the solution slope while only reverting to first-order when necessary.

Vilar \& Abgrall \cite{vilar2024posteriori} introduced a novel concept to further localize changes made to the high-order solution by the limiter. In their approach, a subgrid of FV cell averages is interpolated from any FE of the candidate solution that violate positivity or other prescribed condition. In FV, \textit{a posteriori} limiting is executed on a cell-by-cell basis, identifying to called \textit{troubled cells} whose fluxes need to be revised \cite{clain2011high, loubere2014new, velasco2023spectral}. A robust fallback scheme computes the revision of this intermediate FV solution, after which a the high-degree polynomial of the element is reconstructed from these limited cell averages. This fallback scheme can be the FV MOOD scheme \cite{clain2011high} or any of the other limited schemes discussed here.

The spectral difference (SD) method is another FE approach that can be shown to be equivalent to a quadrature free nodal DG scheme \cite{liu2006spectral, may2011connection}. Velasco \textit{et al.} \cite{velasco2023spectral} took advantage of the fact that the SD method offers a natural interpretation as a nonuniform FV grid, avoiding the need to interpolate between two solution representations. 

For FV schemes, the \textit{a priori} and \textit{a posteriori} limiting paradigms have fundamentally different implementations and trade-offs. Zhang \& Shu's limiter \cite{zhang2011maximum, zhang2010maximum} maintains a strict maximum principle for high-order schemes, but is computationally costly due to the reduced CFL condition and the many nodal reconstructions for the flux quadrature. The cost of the small CFL factor can be mostly mitigated by using an adaptive time-step size \cite{huang2023high}. GPUs, renowned for their ability to execute numerous matrix multiplications as a single vectorized operation, might also potentially mitigate this issue. This is because each nodal reconstruction corresponds to a stencil that can be computed through matrix multiplication.

The cost of the nodal reconstructions is more severe in two dimensions or higher. While the two-dimensional finite volume (FV) flux integral can be reconstructed in various ways, Zhang \& Shu's limiter \cite{zhang2011maximum, zhang2010maximum} use a Gauss-Legendre quadrature. This quadrature method requires a node count per cell face that grows linearly with the polynomial degree. In contrast, transverse flux reconstruction \cite[see e.g.][]{colella2008limiter, mccorquodale2011high, felker2018fourth} uses only one node per cell face, making it a significantly more economical option. We will investigate whether the Zhang \& Shu \cite{zhang2011maximum, zhang2010maximum} \textit{a priori} slope limiter can still be used in an MPP scheme when transverse flux reconstruction is exploited. Previous studies combining an \textit{a priori} limiter with transverse flux reconstruction have resulted in significant violations of the maximum principle \cite{mccorquodale2011high}.

Another limitation of MPP, \textit{a priori} limiters is that, at high-order and after long time-scales, their numerical solutions can be dominated by numerical artifacts and diffusion \cite{kuzmin2022bound}. The high-order schemes in this case actually perform worse at long time-scales than, say, second-order MUSCL-Hancock. Perhaps the trade-off of a strictly MPP scheme is that is too stringent with the permitted slopes as a means to enforces monotonicty. It is unclear how the \textit{a posteriori} limiting compares in this regard.

The \textit{a posteriori} limiters presented here do not require the use of Gauss-Legendre quadrature in two dimensions, making them automatically more cost-effective. However, these limited schemes often violate the maximum principle of the problem at hand, as they lack a guarantee to the contrary. Vilar \& Abgrall \cite{vilar2024posteriori} discovered that extending the convex blending of revised and high-order fluxes to a region of neighbors around each troubled cell reduces the magnitude of these violations. Rueda-Ram\'{i}rez \textit{et al} \cite{rueda2022subcell} showed that a more sophisticated version of this blending can be formulated so as it guarantee the satisfaction of physical bounds on solution variables. Nevertheless, it remains unclear how well a high-order, MPP \textit{a posteriori} limited scheme performs over long time scales and whether they suffer from the same dominance of numerical artifacts as the \textit{a priori} limited schemes.
 
In this study, we develop a methodology to compare \textit{a priori} and \textit{a posteriori} limited explicit finite volume (FV) schemes at arbitrarily high orders. We chose FV as the base scheme over finite element (FE) methods, as modern FE schemes revert to FV when encountering violations, anyways. We solve the linear advection equation in one and two dimensions through a series of numerical experiments. The numerical solutions obtained from the two limiting paradigms are compared with respect to their maximum principle violations and their quality over long time-scales. Furthermore, we evaluate the computational costs of these methods on traditional CPU architecture as well as GPUs.

In Section \ref{sec:numerical_methods}, we give an overview the finite volume method and the various slope-limiting techniques compared in this work. Section \ref{sec:numerical_schemes} conveniently summarizes the \textit{a priori} and \textit{a posteriori} limited schemes studied in our numerical experiments. The results of those experiments are presented in Section \ref{sec:numerical_results}, where we validate our high-order FV and MPP implementations and compare the high-order slope-limited schemes' maximum principle violations, solution quality at long time-scales, and computational cost. Additional commentary regarding the experiments is provided in Section \ref{sec:discussion} and conclusions are drawn in Section \ref{sec:conclusion}.

\section{Numerical methods}
\label{sec:numerical_methods}

We consider the scalar conservation law

\begin{equation}
    \frac{\partial u}{\partial t} + \nabla \cdot \mathbf{f}(u) = 0, \quad u(\mathbf{x}, 0) = u_0(\mathbf{x}) 
    \label{conservation_law}
\end{equation} for $\mathbf{f}, \mathbf{x} \in \mathbb{R}^{d}$, where $u(\mathbf{x},t)$ is a conserved scalar, $\mathbf{f}$ is a flux function and $d$ is the number of dimensions. We begin with the one-dimensional case ($d=1$), and rewrite (\ref{conservation_law}) into the initial value problem

\begin{equation}
    \frac{\partial u}{\partial t} + \frac{\partial f(u)}{\partial x} = 0, \quad u(x, 0)=u_0(x).
    \label{1d_conservation}
\end{equation}

In the finite volume formulation, the spatial domain is partitioned into finite intervals, or cells $I_i=[x_{i-\onehalf}, x_{i+\onehalf}]$ and the cell volume average is defined at time $t^n$ as

\begin{equation}
    \bar{u}^n_i = \frac{1}{h} \int_{x_{i - \onehalf}}^{x_{i + \onehalf}} u(\xi, t^n) d \xi,
    \label{1d_finite_volume}
\end{equation} where $h$ is the length of $I_i$, chosen to be uniform in this work. Combining (\ref{1d_conservation}) with (\ref{1d_finite_volume}) and applying the divergence theorem, the semi-discrete form is written

\begin{equation}
    \frac{d\bar{u}_i^n}{dt} = -\frac{1}{h}(f(u_{i+\onehalf}) - f(u_{i-\onehalf})).
    \label{1d_semidiscrete}
\end{equation} The nodal values $u_{i \pm \onehalf}=P_i(x_{i \pm \onehalf})$ are obtained by reconstructing $u(x,t^n)$ as a polynomial $P_i(x)$ piecewise on each interval $I_i$ and evaluating it at the cell edges. $P_i(x)$ is a conservative interpolation polynomial such that $\int_{I_i}P_i(x)=\bar{u}_i^n$ and is found by taking the derivative of the Lagrange polynomial which interpolates the first moment of $u(x)$ in some kernel of cells containing $I_i$ \cite{harten1997uniformly, shu2009high, zhang2016eno}. The size of the kernel limits the accuracy the polynomial degree of $P_i(x)$ and, consequently, the order of accuracy of the finite volume scheme; $P_i(x)$ is of degree $p$ when it is reconstructed from a kernel of exactly $p + 1$ neighboring cells, so a smooth $u(x)$ is approximated with $u(x) = P_i(x) + \mathcal{O}(h^{p+1})$ on $I_i$.

The evaluation of $P_i(x_{i \pm \onehalf})$ is summarized as stencil operations in Table \ref{tab:conservative_stencils} for various degrees $p$. We select only those stencils that are symmetric so as to simplify our implementation and minimize the number of cells that must be adjusted to enforce boundary conditions. For odd-degree polynomials (which require a kernel of even size), we take the average of the left- and right-biased stencils.

\begin{table}[]
    \centering
    \begin{tabular}{rc}
    \toprule
    $p$ & $u_{i - \frac{1}{2}}^+ = P_i(x_{i - \frac{1}{2}}) = u(x_{i - \frac{1}{2}}) + \mathcal{O}(h^{p+1})$\\
    \midrule
    0 & $\bar{u}_i$ \\
    1 & $\frac{1}{4}\bar{u}_{i-1} + \bar{u}_i - \frac{1}{4}\bar{u}_{i+1}$ \\
    2 & $\frac{1}{3}\bar{u}_{i-1} + \frac{5}{6}\bar{u}_i - \frac{1}{6}\bar{u}_{i+1}$ \\
    3 & $-\frac{1}{24}\bar{u}_{i-2} + \frac{5}{12}\bar{u}_{i-1} + \frac{5}{6}\bar{u}_i - \frac{1}{4}\bar{u}_{i+1} + \frac{1}{24}\bar{u}_{i+2}$ \\
    4 & $-\frac{1}{20}\bar{u}_{i-2} + \frac{9}{20}\bar{u}_{i-1} + \frac{47}{60}\bar{u}_i - \frac{13}{60}\bar{u}_{i+1} + \frac{1}{30}\bar{u}_{i+2}$ \\
    5 & $\frac{1}{120}\bar{u}_{i-3} - \frac{1}{12}\bar{u}_{i-2} + \frac{59}{120}\bar{u}_{i-1} + \frac{47}{60}\bar{u}_i - \frac{31}{120}\bar{u}_{i+1} + \frac{1}{15}\bar{u}_{i+2} - \frac{1}{120}\bar{u}_{i+3}$ \\
    6 & $\frac{1}{105}\bar{u}_{i-3} - \frac{19}{210}\bar{u}_{i-2} + \frac{107}{210}\bar{u}_{i-1} + \frac{319}{420}\bar{u}_i - \frac{101}{420}\bar{u}_{i+1} + \frac{5}{84}\bar{u}_{i+2} - \frac{1}{140}\bar{u}_{i+3}$ \\
    7 & $-\frac{1}{560}\bar{u}_{i-4} + \frac{17}{840}\bar{u}_{i-3} - \frac{97}{840}\bar{u}_{i-2} + \frac{449}{840}\bar{u}_{i-1} + \frac{319}{420}\bar{u}_i - \frac{223}{840}\bar{u}_{i+1} + \frac{71}{840}\bar{u}_{i+2} - \frac{1}{56}\bar{u}_{i+3} + \frac{1}{560}\bar{u}_{i+4}$ \\
    \bottomrule
    \end{tabular}
    \caption{Conservative polynomial interpolation of $u(x)$ at the left cell node $x_{i-\onehalf}$ for degree $p=0,\ldots,7$. $u(x_{i+\onehalf})$ is interpolated by reversing the order of the stencil.}
    \label{tab:conservative_stencils}
\end{table}

To resolve the Riemann problem presented by the nodal values $u_{i+\frac{1}{2}}^-$ and $u_{i+\frac{1}{2}}^+$ interpolated for cells $I_i$ and $I_{i+1}$, respectively, a numerical flux function is written as $\hat{f}_{i \pm \frac{1}{2}} = \mathcal{F}(u_{i + \frac{1}{2}}^-, u_{i + \frac{1}{2}}^+)$. $\mathcal{F}$ must be a Lipschitz continuous function of both arguments, physically consistent with $u$ ($F(u,u)=u$), and monotone, which is to say that it is non-decreasing in its first argument and non-increasing in its second argument \cite{crandall1980monotone}. We use here

\begin{equation}
    \mathcal{F}(u,v) = \frac{1}{2} [ f(u) + f(v) - \max_u(|f'(u)|) (v - u) ],
    \label{GLF}
\end{equation} which satisfies these requirements.

In the case of linear scalar advection with uniform advection velocity $a$, given by

\begin{equation}
    f(u)=au,
    \label{1d_advection}
\end{equation} the Riemann solver given by (\ref{GLF}) is known as the \textit{upwind} solution.

The semi-discrete scheme (\ref{1d_semidiscrete}) is rewritten with numerical fluxes as

\begin{equation}
    \frac{d\bar{u}_i^n}{dt} = \mathcal{L}_p(\bar{u}_i^n) = -\frac{1}{h} \left( \hat{f}_{i+\frac{1}{2}} - \hat{f}_{i-\frac{1}{2}} \right)
    \label{1d_spatial_operator}
\end{equation} and the definition of the order $p+1$ finite volume dynamics, given by the high-order spatial discretization $\mathcal{L}_p$, is complete.

\subsection{Time discretization}

The ordinary differential equation (ODE) in (\ref{1d_spatial_operator}) has the general form

\begin{equation}
    \frac{du}{dt} = D(t, u).
    \label{ode}
\end{equation} The stability of the numerical solution of (\ref{ode}) obtained from Runge-Kutta integration techniques depends on both $D$ as well as the particular Runge-Kutta method used. We will present different Runge-Kutta schemes in the form of Butcher tableaus, defined as follows:

Given times $t^n$ and $t^{n+1} = t^n + \Delta t$, the $s$-stage Runge-Kutta update for $u$ is written as:

\begin{equation}
    u^{n+1} = u^n + \Delta t \sum_{i=1}^s b_i k_i,
\end{equation} where

\begin{align*}
    k_i = D(t^n + c_i \Delta t, u^n + \Delta t \sum_{j=1}^s a_{ij}k_j),
\end{align*} and $b_i$, $c_i$, and $a_{ij}$ are given in the form of Table \ref{tab:RK}.

\begin{table}[]
    \[
    \renewcommand\arraystretch{1.2}
    \begin{array}
    {c|cccc}
    c_1 & a_{11} & a_{12} & \cdots & a_{1s} \\
    c_2 & a_{21} & a_{22} & \cdots & a_{2s} \\
    \vdots & \vdots & \vdots & \ddots & \vdots \\
    c_s & a_{s1} & a_{s2} & \cdots & a_{ss} \\
    \hline
    & b_1 & b_2 & \cdots & b_s
    \end{array}
    \]
    \caption{Generic Butcher tableau for an $s$-stage Runge-Kutta integrator.}
    \label{tab:RK}
\end{table}

The stabilty of Runge-Kutta schemes are assessed with the Dalquist test. Here, we assume linear dynamics $D(u)=\Omega u$, where $\Omega \in \mathbb{C}$ are the eigenvalues of $D$ in (\ref{ode}) \cite{hippolyte2019order, ogunniran2020linear}. The values of $z=\Omega \Delta t$ that result in $|u^{n+1}| \le |u^n|$ give the linear stability region for each integration scheme.

To assess the stability of the spatial discretization in (\ref{1d_spatial_operator}), we analyze the modified wavenumber of $\mathcal{L}_p$ after assuming a harmonic solution $\bar{u}_i^n=e^{ikx}$ \cite{moin2010fundamentals}. This amplification factor is translated to $\Omega \Delta t$-space through the Courant-Friedrichs-Lewy (CFL) factor, given by $C=a \frac{\Delta t}{h}$, where we set $C=1$, for reference.

\begin{figure}
    \centering
    \includegraphics[width=0.9\textwidth]{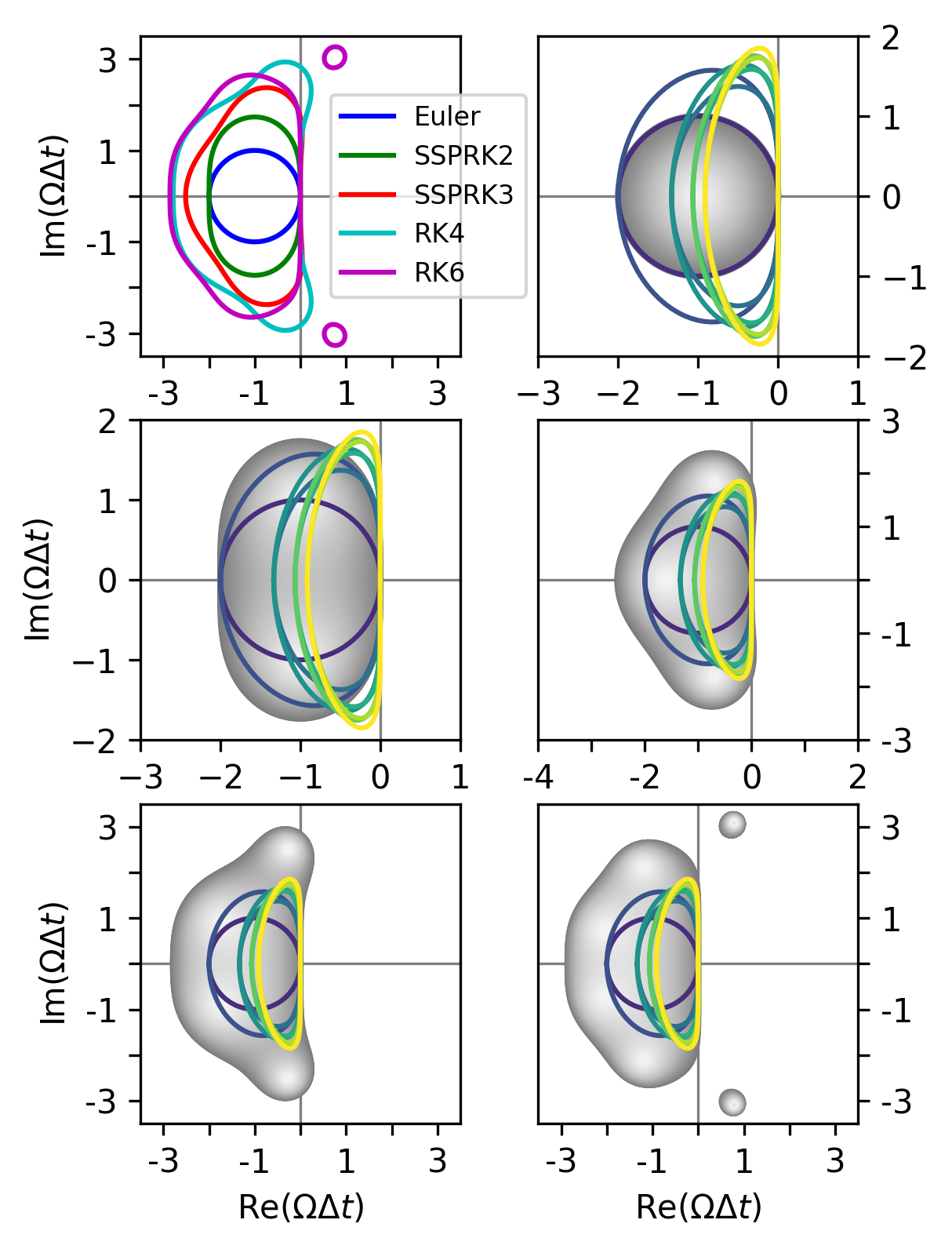}
    \caption{The edges of the stability regions of forward Euler, SSPRK2, SSPRK3, RK4, and RK6 are shown in the top left panel. The subsequent five panels show each Runge-Kutta method overlayed with the eigenvalue tracks of $D=\mathcal{L}_p$ corresponding to finite volume upwinding and a CFL factor $C=1$. The tracks are shown for $p$ from 0 to 7, linearly shaded from purple to yellow.}
    \label{fig:stability}
\end{figure}

The fully-discrete scheme combines (\ref{1d_spatial_operator}) with a Runge-Kutta scheme. If the complete range of modified wavenumbers of $\mathcal{L}_p$ falls entirely within the stability region of a Runge-Kutta scheme, then the fully-discrete scheme is considered stable. For example, the forward Euler method, defined in Table \ref{tab:Euler}, is shown in Figure \ref{fig:stability}. Only the first-order spatial discretization $\mathcal{L}_0$ is contained within the stability region. Any higher-order $\mathcal{L}_p$ combined with the forward Euler method results in numerical instabilities.

The second- and third-order Strong Stability Preserving Runge-Kutta methods (referred to as SSPRK2 and SSPRK3, respectively) are defined in Tables \ref{tab:SSPRK2} and \ref{tab:SSPRK3}, where the Strong Stability Preserving property ensures a strict maximum principle is maintained in the numerical solution \cite{zhang2011maximum, zhang2010maximum, sun2019strong, hadjimichael2013strong, gottlieb2009high}. The maximum principle of (\ref{conservation_law}) is given by $M=\max_\mathbf{x} u_0(\mathbf{x})$ and $m=\min_\mathbf{x} u_0(\mathbf{x})$ if $u(\mathbf{x}, t) \in [m, M]$ for all $\mathbf{x}$ and $t$ \cite{zhang2011maximum, zhang2010maximum}. The first-order upwinding scheme $\mathcal{L}_0$ paired with the forward Euler method is termed maximum-principle-preserving (MPP) as it guarantees $\bar{u}_i^n \in [m,M]$ throughout the computational domain \cite{zhang2011maximum, zhang2010maximum}. Furthermore, any time integration method that is equivalent to a convex combination of forward Euler steps (like SSPRK2 and SSPRK3) will also be MPP for $\mathcal{L}_0$ \cite{zhang2011maximum, zhang2010maximum}.

Due to the Godunov order barrier theorem, any $\mathcal{L}_p$ with $p>0$ will not be MPP when used with any linear numerical ODE method. This issue is addressed later with additional numerical techniques.

Figure \ref{fig:stability} shows that SSPRK2 is stable for $\mathcal{L}_0$, $\mathcal{L}_1$. It is only apparent from closely zooming into $\text{Re} (\Omega \Delta t)=0$ that the eigenvalue track of $\mathcal{L}_2$ falls slightly outside the stability region of SSPRK2. On the other hand, SSPRK3 and all higher-order Runge-Kutta methods provided, remain stable for $\mathcal{L}_0, $ through $\mathcal{L}_7$. 

The classic fourth-order Runge-Kutta method (RK4), as defined in Table \ref{tab:RK4}, is not a convex combination of forward Euler steps. Despite this, it has been demonstrated to be quasi-SSP \cite{sanz2010positivity}. We will further investigate how RK4 performs when used with MPP and approximately-MPP spatial discretization methods that are outlined in subsequent subsections.

The highest time integration method used in this work is sixth-order Runge-Kutta (RK6), given in Table \ref{tab:RK6} \cite{luther1968explicit}. While it is not anticipated to maintain a strict maximum principle when used in conjunction with any spatial discretization, it will be useful for conducting tests involving very high-order schemes.

\begin{table}[]
    \[
    \renewcommand\arraystretch{1.2}
    \begin{array}
    {c|c}
    0 & 0 \\
    \hline
    & 1
    \end{array}
    \]
    \caption{Butcher tableau for explicit, first-order Euler integration.}
    \label{tab:Euler}
\end{table}

\begin{table}[]
    \[
    \renewcommand\arraystretch{1.2}
    \begin{array}
    {c|cc}
    0 & 0 & 0 \\
    1 & 1 & 0\\
    \hline
    & \frac{1}{2} &\frac{1}{2}
    \end{array}
    \]
    \caption{Butcher tableau for explicit, second-order Strong Stability Preserving Runge-Kutta integration (SSPRK2).}
    \label{tab:SSPRK2}
\end{table}

\begin{table}[]
    \[
    \renewcommand\arraystretch{1.2}
    \begin{array}
    {c|ccc}
    0 & 0 & 0 & 0 \\
    1 & 1 & 0 & 0 \\
    \frac{1}{2} & \frac{1}{4} & \frac{1}{4} & 0 \\
    \hline
    & \frac{1}{6} &\frac{1}{6} & \frac{2}{3}
    \end{array}
    \]
    \caption{Butcher tableau for explicit, third-order Strong Stability Preserving Runge-Kutta integration (SSPRK3).}
    \label{tab:SSPRK3}
\end{table}

\begin{table}[]
    \[
    \renewcommand\arraystretch{1.2}
    \begin{array}
    {c|cccc}
    0 & 0 & 0 & 0 & 0 \\
    \frac{1}{2} & \frac{1}{2} & 0 & 0 & 0 \\
    \frac{1}{2} & 0 & \frac{1}{2} & 0 & 0 \\
    1 & 0 & 0 & 1 & 0 \\
    \hline
    & \frac{1}{6} &\frac{1}{3} &\frac{1}{3} &\frac{1}{6} 
    \end{array}
    \]
    \caption{Butcher tableau for explicit, fourth-order Runge-Kutta integration (RK4).}
    \label{tab:RK4}
\end{table}

\begin{table}[]
    \[
    \renewcommand\arraystretch{1.2}
    \begin{array}
    {c|ccccccc}
    0 & 0 & 0 & 0 & 0 & 0 & 0 & 0 \\
    
    1 & 1 & 0 & 0 & 0 & 0 & 0 & 0 \\
    
    \frac{1}{2} & \frac{3}{8} & \frac{1}{8} & 0 & 0 & 0 & 0 & 0 \\
    
    \frac{2}{3} & \frac{8}{27} & \frac{2}{27} & \frac{8}{27} & 0 & 0 & 0 & 0 \\
    
    \frac{7 - \sqrt{21}}{14} & \frac{3(3\sqrt{21}-7)}{392} & \frac{\sqrt{21} - 7}{49} & \frac{6(7 - \sqrt{21})}{49} & \frac{-3(21 - \sqrt{21})}{392} & 0 & 0 & 0 \\
    
    \frac{7 + \sqrt{21}}{14} & \frac{-231 - 51\sqrt{21}}{392} & \frac{-7 - \sqrt{21}}{49} & \frac{-8\sqrt{21}}{49} & \frac{3(21+121\sqrt{21})}{1960} & \frac{49(6+\sqrt{21})}{245} & 0 & 0 \\
    
    1 & \frac{22+7\sqrt{21}}{12} & \frac{2}{3} & \frac{2(7\sqrt{21} - 5)}{9} & \frac{-63(3\sqrt{21} - 2)}{180} & \frac{-7(49+9\sqrt{21})}{90} & \frac{7(7-\sqrt{21})}{18} & 0 \\
    
    \hline
    & \frac{1}{20} & 0 & \frac{16}{45} & 0 & \frac{49}{180} & \frac{49}{180} & \frac{1}{20} 
    \end{array}
    \]
    \caption{Butcher tableau for explicit, sixth-order Runge-Kutta integration (RK6).}
    \label{tab:RK6}
\end{table}

\subsection{\textit{A priori} slope limiting}

We introduce the concept of a slope limiter, which is a modification to the basic finite volume method described in (\ref{1d_spatial_operator}), enabling higher-order ($p>0$) solvers to achieve MPP or approximately-MPP behavior. In the presence of discontinuities, higher-order interpolation polynomials can exhibit unphysical oscillations. These oscillations disrupt the local monotonicity of the numerical solution and may cause violations of the maximum principle. Slope limiters aim to mitigate these issues by locally reducing the polynomial degree of the spatial discretization where oscillations are detected, thereby ensuring a smoother gradient at discontinuities. Importantly, slope limiters are designed to maintain the conservative property of the base finite volume scheme.

In this work, we focus on the distinction between \textit{a priori} and \textit{a posteriori} slope limiters. \textit{A priori} slope limiters detect oscillations from only the data at $t^n$, which is used to compute the slope-limited update at $t^{n+1}$. On the other hand, \textit{a posteriori} slope limiters involve computing a candidate solution at $t^{n+1}$ and computing revisions of this tentative solution after the facts.

Zhang \& Shu \cite{zhang2011maximum, zhang2010maximum} developed an \textit{a priori} finite volume slope limiter that is MPP at arbitrary order. Their method relies on a parameter $\theta_i \in [0,1]$, which blends high-order and first-order interpolations of the nodal values within each cell, resulting in the following modified polynomial

\begin{equation}
    \Tilde{P}_i(x) = \theta_i \left( P_i(x) - \bar{u}_i \right) + \bar{u}_i, \quad \theta_i = \min \left( \left| \frac{M_i - \bar{u}_i}{M'_i - \bar{u}_i} \right|, \left| \frac{m_i - \bar{u}_i}{m'_i - \bar{u}_i} \right|, 1 \right).
 	\label{apriori_slope_limiter}
\end{equation} Here, $M_i=\max (\bar{u}_{i-1}, \bar{u}_i, \bar{u}_{i+1})$ and $m_i=\min (\bar{u}_{i-1}, \bar{u}_i, \bar{u}_{i+1})$ are the maximum and minimum of each cell and its adjacent neighbors. $M'=\max_{x \in S_i}P_i(x)$ and $m'=\min_{x \in S_i}P_i(x)$ denote the maximum and minimum of $P_i(x)$ evaluated on the set of the $L$-point Gauss-Lobatto quadrature points $S_i = \{ x_{i-\onehalf} = x_i^1, x_i^2, \ldots, x_i^{L-1}, x_i^{L} = x_{i+\onehalf} \}$.

It is evident from (\ref{apriori_slope_limiter}) that the high-order nodal values are limited when $\theta_i < 1$. This happens when the difference between the unlimited polynomial interpolation of $u(t,x)$ over $I_i$ and the first-order interpolation $\bar{u}_i$ is greater than the difference between the local maximum principle, represented by $M_i$ and $m_i$, and $\bar{u}_i$ --- precisely what occurs in the presence of a spurious oscillation.

The Gauss-Lobatto quadrature rule for computing an integral is exact for polynomials of degree $2L - 3$ or less, so $L$ is chosen to be the smallest integer such that $2L - 3 \ge p$. Since $p>1$ requires nodes other than $u(x_{i \pm \onehalf})$, we need additional stencils to the ones provided in Table \ref{tab:conservative_stencils} in order to compute $\theta_i$. We will address this later with a slight modification of the original method.

Zhang \& Shu's \textit{a priori} slope-limited spatial discretization is guaranteed to satisfy the maximum principle if it is solved with an SSP Runge-Kutta method and also if it satisfies a reduced CFL condition \cite{zhang2011maximum, zhang2010maximum}. Let $w_\alpha$ be the weights of the $L$-point Gauss-Lobatto quadrature rule such that $\sum_{\alpha=1}^L w_\alpha = 1$. We write this reduced CFL condition as

\begin{equation}
    C < C_\text{MPP} = \min[w_0, \ldots, w_L].
    \label{reduced_courant}
\end{equation} Values of $C_\text{MPP}$ are provided in Table \ref{tab:reduced_courant}.

\begin{table}[]
    \centering
    \begin{tabular}{ccc}
        \toprule
        $p$ & $C_\text{MPP}$ & Quadrature points on $[-\frac{1}{2}, \frac{1}{2}]$\\
        \midrule
        0, 1 & $\frac{1}{2}$ & $\{ \pm \frac{1}{2} \}$ \\
        2, 3 & $\frac{1}{6}$ & $\{ 0, \pm \frac{1}{2} \}$ \\
        4, 5 & $\frac{1}{12}$ & $\{ \pm \frac{1}{\sqrt{20}}, \pm \frac{1}{2} \}$ \\
        6, 7 & $\frac{1}{20}$ & $\{ 0, \pm \frac{\sqrt{21}}{14}, \pm \frac{1}{2} \}$ \\
        \bottomrule
    \end{tabular}
    \caption{The maximum CFL factor $C_\text{MPP}$ and the corresponding $L$ Gauss-Lobatto quadrature points used in the Zhang \& Shu slope limiter for polynomial degree $p=0,\ldots,7$, where $L$ is chosen as the smallest value such that $2L-3 \ge p$.}
    \label{tab:reduced_courant}
\end{table}

The reduced $C$ (resulting in a reduced $\Delta t$) and the increase in the required number of interpolated values make the \textit{a priori} limiting scheme significantly more costly at higher order. To mitigate this issue, we propose two modifications. 

\begin{itemize}
\item First, we suggest simplifying the set of points used to compute $\theta_i$ for $p>3$. Instead of computing $M'_i$ and $m'_i$ by finding the minimum and maximum of $P_i(x)$ evaluated at the $L>4$ Gauss-Lobatto quadrature points $S_i$, we use only three points $\{ x_{i-\onehalf}, x_i, x_{i+\onehalf} \}$, where $x_i$ is the centroid of the cell. This reduced set of points is referred to as the \textit{centroid} set.
\item The second modification involves adopting an adaptive time-step size. The reduced CFL condition (\ref{reduced_courant}) requires significantly smaller $\Delta t$ at higher order. However, in practice, the small time-step size is only necessary for the initial steps when the numerical solution features sharp continuities. As the numerical solution becomes progressively smoother, larger time-step sizes can be used while maintaining a strict maximum principle. We follow the adaptive time-stepping technique proposed by \cite{huang2023high}, where $\bar{u}^{n+1}$ is recomputed if it violates the maximum principle with $\frac{1}{2} \Delta t$. This process is repeated until the solution is MPP, or once $\Delta t$ satisfies (\ref{reduced_courant}), at which point the step will be MPP, regardless. In this modification, we initially attempt $C=0.8$ at each step.
\end{itemize}

\subsection{\textit{A posteriori} slope limiting}

Another approach to the problem of slope limiting is the \textit{a posteriori} limiter, first implemented by Clain \textit{et al} as Multi-dimensional Optimal Order Detection (MOOD) \cite{clain2011high}. Here, a candidate solution at $t^{n+1}$ is computed and screened for \textit{troubled cells}. These cells are identified based on violations of the maximum principle or other predefined criteria. Subsequently, the solutions at troubled cells are locally recomputed using a lower-order scheme, generating a new candidate solution. This new solution is then reevaluated for troubled cells, and the process continues until either no troubled cells remain or there are no lower-order schemes with which to revise them.

An alternative approach, proposed by \cite{velasco2023spectral} for the Spectral Difference (SD) finite element framework, jumps straight to a robust, second-order fallback scheme to recompute troubled cells. Since there is only one revision step in this implementation, it skips the cumbersome recomputing of solutions with one less polynomial degree at a time that is required by MOOD. The fallback scheme used in \cite{velasco2023spectral} is the well-established second-order MUSCL-Hancock reference scheme, which is strictly MPP. The methodology of \cite{velasco2023spectral} functions like MOOD with only two options for the so-called ``cascade of schemes'': $p>1$ and $p=1$.

Regardless of the implementation details of the \textit{a posteriori} limiter, the underlying concept is akin to that of \textit{a priori} slope limiting: when an interpolation polynomial exhibits unphysical oscillations, resulting in over or undershoots beyond a permissible range, the local degree of the interpolation polynomial should be reduced in such problematic regions.

In this work, we take inspiration from \cite{velasco2023spectral} and use MUSCL-type schemes to revise our high-order candidate solutions at each step if needed.

Let $\bar{u}^*_i$ denote the candidate solution of cell $i$ for time $t^{n+1}$. A cell is flagged as troubled if it does not satisfy the condition for \textit{numerical admissibility detection} (NAD):

\begin{equation}
    m_i - \epsilon (M-m) \le \bar{u}^*_i \le M_i + \epsilon (M-m),
    \label{NAD}
\end{equation} where $M_i$ and $m_i$ are defined the local maximum principle from the previous subsection and $\epsilon$ is a small tolerance, chosen to be $1\times10^{-5}$ \cite{velasco2023spectral}. $\epsilon$ is useful in regions where the numerical solution is uniform, since the local maximum principle $M_i$ and $m_i$ will otherwise be sensitive to the slight oscillations, flagging more cells as troubled than needed.

\subsubsection{MUSCL Fallback scheme}

The nodal interpolations given by a second-order MUSCL scheme can be written

\begin{equation}
    \Tilde{u}_{i-\onehalf}^+ = \bar{u}_i - \onehalf \Tilde{S}_i, \quad \Tilde{u}_{i+\onehalf}^- = \bar{u}_i + \onehalf \Tilde{S}_i,
\end{equation} where $\Tilde{S}_i$ is a modification of the local slope multiplied by $h$:

\begin{equation}
    \Tilde{S}_i=\text{lim}(\bar{u}_{i-1}, \bar{u}_i, \bar{u}_{i+1}) \frac{\bar{u}_{i+1}-\bar{u}_{i-1}}{2}.
\end{equation} The slope limiter satisfies

\begin{equation}
    0 \le \text{lim}(\bar{u}_{i-1}, \bar{u}_i, \bar{u}_{i+1}) \le 1,
\end{equation} where a value of 1 gives the degree $p=1$ nodal interpolation and 0 gives $p=0$ (see Table \ref{tab:conservative_stencils}), so it functions exactly like $\theta$, the high-order \textit{a priori} slope limiter.

Let $k$ be the index of a troubled cell. We compute the fluxes given by the fallback scheme

\begin{equation}
    \begin{split}
        \hat{f}_{k-\onehalf}^M = \mathcal{F}( \Tilde{u}_{k-\onehalf}^-, \Tilde{u}_{k-\onehalf}^+ ) \\
        \hat{f}_{k+\onehalf}^M = \mathcal{F}( \Tilde{u}_{k+\onehalf}^-, \Tilde{u}_{k+\onehalf}^+ )
    \end{split}
\end{equation} and reassign high-order fluxes from (\ref{1d_spatial_operator}) with $\hat{f}_{k \pm \onehalf} \leftarrow \hat{f}^M_{k \pm \onehalf}$.

The fallback scheme is only second-order when coupled with a second-order or higher Runge-Kutta method. In this setup, we search for troubled cells and compute their revised fluxes at each integration stage. This requires computing a candidate solution for each stage, which is achieved by performing an Euler step with a time-step size given for each of the $s$ stages as $c_1, \ldots, c_s$ in chosen Runge-Kutta method's Butcher Tableau.

%\subsubsection{MUSCL-Hancock Fallback scheme}
%
%Alternatively, in the spirit of MUSCL-Hancock, we can use at each revision step a predictor-corrector scheme. The predictor step is found with the half time-step size update
%
%\begin{equation}
%    \bar{u}_i^{n+\onehalf} = \bar{u}_i^n - \frac{a \Delta t}{2h} S_i,
%\end{equation} then the slope-limited nodal interpolations are found with the correction
%
%\begin{equation}
%    \Tilde{u}_{i-\onehalf}^+ = \bar{u}_i^{n+\onehalf} - \onehalf \Tilde{S}_i, \quad \Tilde{u}_{i+\onehalf}^- = \bar{u}_i^{n+\onehalf} + \onehalf \Tilde{S}_i.
%\end{equation} Again, here we use $\Delta t$ given by the corresponding stage of the Butcher Tableau.

The left and right differences of cell $i$ are defined:

\begin{equation}
    S_i^L = \bar{u}_i - \bar{u}_{i-1}, \quad S_i^R = \bar{u}_{i+1} - \bar{u}_i.
\end{equation} Then, the \textit{minmod}-limited difference is given by

\begin{equation}
    \Tilde{S}_i = \text{sgn}(S_i) \cdot \min(|S_i^L|, |S_i^R|).
\end{equation} Note that the MUSCL and MUSCL-Hancock schemes are strictly MPP when using this limiter.

The central differences of cell $i$ is defined:

\begin{equation}
    S_i^C = \frac{1}{2}(S_i^L + S_i^R).
\end{equation} Then, the \textit{moncen}-limited difference is given by

\begin{equation}
    \Tilde{S}_i = \text{sgn}(S_i) \cdot \min(|2S_i^L|, |S_i^C|, |2S_i^R|).
\end{equation} The \textit{moncen} limiter is strictly MPP like \textit{minmod} with the added advantage of being significantly less diffusive.

\subsubsection{Blended flux correction}

\begin{figure}
    \centering
    \includegraphics[width=\textwidth]{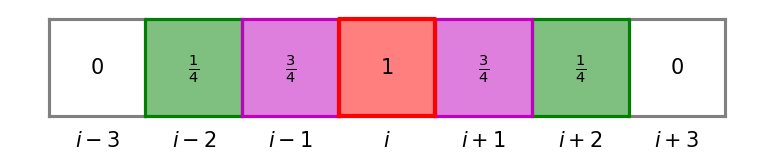}
    \caption{The value of $\phi$ used for convex blending of corrected fluxes, shown across a region comprising seven cells, with one troubled cell highlighted in red.}
    \label{fig:blending1d}
\end{figure}

Recent work by Vilar \& Abgrall \cite{vilar2024posteriori} found that blending the high-order and fallback fluxes of the cells neighboring one that is troubled improves the preservation of the maximum principle. Instead of assigning a high-order or fallback flux at each cell interface, we assign a convex combination of the two:

\begin{equation}
    \begin{split}
        \hat{f}_{i - \onehalf} \leftarrow \max(\phi_{i-1}, \phi_i) (\hat{f}^M_{i-\onehalf} - \hat{f}_{i-\onehalf}) + \hat{f}_{i-\onehalf} \\
        \hat{f}_{i + \onehalf} \leftarrow \max(\phi_i, \phi_{i+1}) (\hat{f}^M_{i+\onehalf} - \hat{f}_{i+\onehalf}) + \hat{f}_{i+\onehalf},
    \end{split}
    \label{convexblending}
\end{equation} where $\phi_i \in [0,1]$ is computed at each cell with the so-called \textit{naive} procedure outline by \cite{vilar2024posteriori}. The procedure is the following. If cell $i$ is troubled, it is given $\phi_i=1$. If it is not troubled and immediately adjacent to a troubled cell, it is given $\phi_i=\frac{3}{4}$. If a cell is not troubled and located two cell positions away from the nearest troubled cell, it is given $\phi_i=\frac{1}{4}$. Finally, if a cell is not troubled and three or more cell positions away from the nearest troubled cell, it is given $\phi_i=0$. An example is shown in Figure \ref{fig:blending1d}.

Note that the non-blended flux correction is recovered by setting $\phi_i=1$ when cell $i$ is troubled and $\phi_i=0$ otherwise.

\subsection{Smooth extrema detection}

The well-known drawback of \textit{a priori} and \textit{a posteriori} slope limiters alike is that they can reduce the accuracy of the numerical solution in smooth, well-behaved regions where they are not needed. To combat this issue, we implement \textit{smooth extrema detection} (SED), which modifies the slope limiter such that it deactivates at cells which are flagged as smooth extrema.

The derivative of the solution is evaluated using a central finite differences approximation:

\begin{equation}
	\bar{u}'_i = \bar{u}_{i+1} - \bar{u}_{i-1}.
\end{equation} Central, left, and right differences are taken from these first differences, giving 

\begin{equation}
	S_i^C = \frac{\bar{u}'_{i+1} - \bar{u}'_{i-1}}{2}, \quad
	S_i^L = \bar{u}'_i - \bar{u}'_{i-1}, \quad
	S_i^R = \bar{u}'_{i+1} - \bar{u}'_i.
\end{equation} With these second derivative terms, we define the left and right smoothness indicators

\begin{equation}
    \alpha^{L,R}_i =
    \begin{cases}
        \min \left( 1, \min \left( 2 S^{L,R}_i, 0 \right) \frac{1}{S^C_i} \right) & S^C_i < 0 \\
        1 & S^C_i = 0 \\
        \min \left( 1, \max \left( 2 S^{L, R}_i, 0 \right) \frac{1}{S^C_i} \right) & S^C_i > 0.
    \end{cases}
\end{equation} The overall smoothness indicator of cell is the minimum of the two:

\begin{equation}
    \alpha_i = \min ( \alpha^L_i, \alpha^R_i ).
    \label{smooth_extrema_detection}
\end{equation} Then, we flag cell $i$ as a smooth extremum if and only if

\begin{equation}
    \min (\alpha_{i-1}, \alpha_i, \alpha_{i+1}) = 1.
\end{equation}

This implementation can create small violations of the maximum principle if it permits undershoots or overshoots in smooth regions. In such cases, we ignore the suggestion of the SED algorithm and use the limited slope. Since we want to avoid the excessive sensitivity of our slope limiting in regions where the solution is uniform, we introduce an additional small tolerance $\epsilon_m$, chosen to be $1\times10^{-10}$. The SED procedures for \textit{a priori} and \textit{a posteriori} slope limiting are given in Algorithms \ref{alg:SED_apriori} and \ref{alg:SED_aposteriori}, respectively.

\begin{algorithm}
\caption{Smooth extrema detection for the \textit{a priori} slope limiter.}\label{alg:SED_apriori}
\begin{algorithmic}
    \State Compute the \textit{a priori} slope limiting terms $M'_i$, $m'_i$, $\theta_i$ at all cells.
    \State Compute the smoothness indicator $\alpha_i$ at all cells.
    \If{$M'_i > M + \epsilon_m \text{ or } m'_i < m - \epsilon_m$}
        \State Do not modify $\theta_i$.
    \Else
        \If{$\min (\alpha_{i-1}, \alpha_i, \alpha_{i+1}) = 1$}
            \State $\theta_i \leftarrow 1$
        \Else
            \State Do not modify $\theta_i$.
        \EndIf
    \EndIf
\end{algorithmic}
\end{algorithm}

\begin{algorithm}
\caption{Smooth extrema detection for the \textit{a posteriori} slope limiter.}\label{alg:SED_aposteriori}
\begin{algorithmic}
    \State Compute the candidate solution $\bar{u}_i^*$ at all cells. 
    \State Compute the smoothness indicator $\alpha_i$ at all cells from this candidate solution.
    \State Flag all cells as \textit{troubled} or \textit{not troubled} based on NAD.
    \If{$\bar{u}_i^* > M + \epsilon_m \text{ or } \bar{u}_i^* < m - \epsilon_m$}
        \State Flag cell $i$ as troubled.
    \Else
        \If{$\min (\alpha_{i-1}, \alpha_i, \alpha_{i+1}) = 1$}
            \State Flag cell $i$ as not troubled.
        \Else
            \State Do not modify not modify the flag of cell $i$.
        \EndIf
    \EndIf
\end{algorithmic}
\end{algorithm}

\subsection{Modifications of the slope-limited, finite volume schemes for two dimensions}

The two-dimensional ($d=2$) form of (\ref{conservation_law}) writes

\begin{equation}
    \frac{\partial u}{\partial t} + \frac{\partial f(u)}{\partial x} + \frac{\partial g(u)}{\partial y} = 0,
\end{equation} where $f$ and $g$ are one-dimensional flux functions in the $x$- and $y$-directions, respectively. 

For the finite volume formulation, we partition the domain into intervals $I_{i,j} = [x_{i-\frac{1}{2}}, x_{i+\frac{1}{2}}] \cup [y_{j-\frac{1}{2}}, y_{j+\frac{1}{2}}]$. Then, the finite volume cell average is defined as

\begin{equation}
    \overline{\overline{u}}_{i, j}^n = \frac{1}{h^2} \int_{y_{j-\frac{1}{2}}}^{y_{j+\frac{1}{2}}} \int_{x_{i-\frac{1}{2}}}^{x_{i+\frac{1}{2}}} u(\xi, \eta, t^n) d\xi d\eta,
\end{equation} where the interval side lengths, chosen to be equal and uniform in this work, are given by $h$.

The semi-discrete scheme in two-dimensions is written as

\begin{equation}
    \begin{split}
        \frac{d\overline{{\overline{{u}}}}_{i,j}^n}{dt} = -\frac{1}{h^2} \int_{y_{j-\onehalf}}^{y_{j+\onehalf}} \hat{f}[u_{i+\onehalf,j}^-(\eta), u_{i+\onehalf,j}^+(\eta)] - \hat{f}[u_{i-\onehalf,j}^-(\eta), u_{i-\onehalf,j}^+(\eta)] d\eta \\ - \frac{1}{h^2} \int_{x_{i-\onehalf}}^{x_{i+\onehalf}} \hat{g}[u_{i,j+\onehalf}^-(\xi), u_{i,j+\onehalf}^+(\xi)] - \hat{g}[u_{i,j-\onehalf}^-(\xi), u_{i,j-\onehalf}^+(\xi)] d\xi,
    \end{split}
\label{semidiscrete2d}
\end{equation} where $\hat{f}(u,v)$ and $\hat{g}(u,v)$ are the numerical fluxes. The traces of $u(x,y,t^n)$ along each cell face are given by the cell-wise conservative polynomial reconstruction $P_{i,j}(x,y)$. Namely,

\begin{equation}
    \begin{split}
        u_{i-\onehalf,j}^+(y) = P_{i,j}(x_{i-\onehalf}, y) \\
        u_{i+\onehalf,j}^-(y) = P_{i,j}(x_{i+\onehalf}, y) \\
        u_{i,j-\onehalf}^+(x) = P_{i,j}(x, y_{i-\onehalf}) \\
        u_{i,j+\onehalf}^-(x) = P_{i,j}(x, y_{i+\onehalf}). \\
    \end{split}
\end{equation}

\subsubsection{Gauss-Legendre quadrature}

Each numerical flux in (\ref{semidiscrete2d}) appears as an integral along its respective cell face. We present two options with which to evaluate the integral, the first being the Gauss-Legendre quadrature rule, summarized in Table \ref{tab:gauss}.

Along each cell face, nodal values of $u(x,y)$ corresponding to the positions of the Gauss-Legendre quadrature points are interpolated by the conservative polynomial $P_{i,j}(x,y)$. This interpolation process is equivalent to a stencil operation; While the stencils in Table \ref{tab:conservative_stencils} are presented for the one-dimensional finite volume method, they can still be applied in two-dimensions: first to reconstruct line averages from the two-dimensional cell volume averages and then to interpolate nodal values from those line averages.

Stencil operations corresponding to the reconstruction at the Gauss-Legendre quadrature points are not provided in Table \ref{tab:conservative_stencils}, but can be derived via the same procedure.

The point-wise numerical flux is then computed at each cell face, with each point requiring its own Riemann solution (see Figures \ref{fig:gaussgausslobato} and \ref{fig:gausscentroid}). For example, suppose we interpolate nodal values of $u(x,y)$ along the $y=y_{j + \onehalf}$ face of cell $i,j$ corresponding to the positions of the $K$-point Gauss-Legendre quadrature $\hat{S}_i = \{ \hat{x}_i^\beta : \beta = 1, \ldots, K \}$. Then we must compute $\{ \hat{g}(u_{i,j+\onehalf}^{\beta,-}, u_{i,j+\onehalf}^{\beta, +}): \beta = 1, \ldots, K \}$ where $u_{i,j+\onehalf}^{\beta,-} =  P_{i,j}(\hat{x}_i^\beta, y_{j+\onehalf})$ and $u_{i,j+\onehalf}^{\beta, +} = P_{i,j+1}(\hat{x}_i^\beta, y_{j+\onehalf})$. Finally, the integral along the cell face is given by

\begin{equation}
    \int_{x_{i-\onehalf}}^{x_{i+\onehalf}} \hat{g}[ u_{i,j+\onehalf}^-(\xi) , u_{i,j+\onehalf}^+(\xi)] d\xi = \sum_{\beta=1}^K \hat{w}_\beta \cdot \hat{g}[ u_{i,j+\onehalf}^{\beta, -}, u_{i,j+\onehalf}^{\beta, +} ]
\end{equation}, where $\hat{w_\beta}$ are the $K$-point Gauss-Legendre quadrature points such that $\sum_{\beta=1}^K \hat{w}_\beta = 1$.

The $K$-point Gauss-Legendre quadrature rule is exact for polynomials of degree $2K-1$ or less, so $K$ is chosen as the smallest integer such that $2K-1 \ge p$.

\subsubsection{Transverse flux reconstruction}

Alternatively, the integral of the flux along cell faces can be computed using only one node from adjacent cell faces. In this method, we interpolate the nodal values of $u(x,y)$ at the midpoint of each cell face from the reconstructed conservative polynomial $P_{i,j}(x,y)$ (see Figure \ref{fig:transversecentroid}). Stencil operations corresponding to this midpoint interpolation are given in Table \ref{tab:transverse_midpoint}. For instance, let's take the nodal value $u(\hat{x}_i, y_{j+\onehalf})$, which is defined at the midpoint of the $y_{j+\onehalf}$ face. The single point-wise flux at each $y_{j+\onehalf}$ face is then given by the Riemann problem $\hat{g}[P_{i,j}(\hat{x}_i, y_{j+\onehalf}), P_{i,j+1}(\hat{x}_i, y_{j+\onehalf})]$. The face integral of this numerical flux can be evaluated by constructing a Lagrange interpolation polynomial at the neighboring point-wise fluxes, all uniformly spaced by a distance of $h$. Since the integral of the flux along the cell face depends on the values of cell faces located transverse to itself, this method is referred to as \textit{transverse flux reconstruction}. Stencil operations corresponding to this integral are given in Table \ref{tab:transverse_integral}. The stencil length is always chosen such that it is exact for polynomials of degree $p$.

\begin{table}[]
    \centering
    \begin{tabular}{ccc}
        \toprule
        $p$ & Quadrature points on $[0, \frac{1}{2}]$ & Quadrature weights \\
        \midrule
        0, 1 & $\{ 0 \}$ & $\{ 1 \}$ \\
        2, 3 & $\{\frac{1}{2 \sqrt{3}}\}$ & $\{ \frac{1}{2} \}$ \\
        4, 5 & $\{ 0, \frac{\sqrt{3}}{2 \sqrt{5}} \}$ & $\{ \frac{4}{9}, \frac{5}{18} \}$ \\
        6, 7 & $ \left\{ \sqrt{\frac{3}{7} - \frac{2}{7} \sqrt{\frac{6}{5}}}, \sqrt{\frac{3}{7} + \frac{2}{7} \sqrt{\frac{6}{5}}} \right\} $ & $\{ \frac{18 + \sqrt{30}}{36}, \frac{18 - \sqrt{30}}{36} \}$ \\
        \bottomrule
    \end{tabular}
    \caption{Gauss-Legendre quadrature on $[-\frac{1}{2}, \frac{1}{2}]$ at varying polynomial degrees $p$. Only the quadrature points and weights on $[0, \frac{1}{2}]$ are reported due to their symmetry about $0$.}
    \label{tab:gauss}
\end{table}

\begin{table}[]
    \centering
    \begin{tabular}{cc}
    \toprule
    $p$ & $u_i = P_i(\bar{x}_i) = u(\bar{x}_i) + \mathcal{O}(h^{p+1})$\\
    \midrule
    0, 1 & $\bar{u}_i$ \\
    2, 3 & $-\frac{1}{24}\bar{u}_{i-1} + \frac{13}{12}\bar{u}_i - \frac{1}{24}\bar{u}_{i+1}$ \\
    4, 5 & $\frac{3}{640}\bar{u}_{i-2} - \frac{29}{480}\bar{u}_{i-1} + \frac{1067}{960}\bar{u}_i - \frac{29}{480}\bar{u}_{i+1} + \frac{3}{640}\bar{u}_{i+2}$ \\
    6, 7 & $-\frac{5}{7168}\bar{u}_{i-3} + \frac{159}{17920}\bar{u}_{i-2} - \frac{7621}{107520}\bar{u}_{i-1} + \frac{30251}{26880}\bar{u}_i - \frac{7621}{107520}\bar{u}_{i+1} + \frac{159}{17920}\bar{u}_{i+2} - \frac{5}{7168}\bar{u}_{i+3}$ \\
    \bottomrule
    \end{tabular}
    \caption{Conservative polynomial interpolation of $u(x)$ at a cell midpoint $\bar{x}_i$ for degree $p=0,\ldots,7$.}
    \label{tab:transverse_midpoint}
\end{table}

\begin{table}[]
    \centering
    \begin{tabular}{cc}
    \toprule
    $p$ & $ \int_{x_{i-\onehalf}}^{x_{i+\onehalf}} u(\xi) d\xi = \int_{x_{i-\onehalf}}^{x_{i+\onehalf}} Q_i(\xi) d\xi + \mathcal{O}(h^{p + 1})$ \\
    \midrule
    0, 1 & $u_i$ \\
    2, 3 & $\frac{1}{24}u_{i-1} + \frac{11}{12}u_i + \frac{1}{24}u_{i+1}$ \\
    4, 5 & $-\frac{17}{5760}u_{i-2} + \frac{77}{1440}u_{i-1} + \frac{863}{960}u_i + \frac{77}{1440}u_{i+1} - \frac{17}{5760}u_{i+2}$ \\
    6, 7 & $\frac{367}{967680}u_{i-3} - \frac{281}{53760}u_{i-2} + \frac{6361}{107520}u_{i-1} + \frac{215641}{241920}u_i + \frac{6361}{107520}u_{i+1} - \frac{281}{53760}u_{i+2} + \frac{367}{967680}u_{i+3}$ \\
    \bottomrule
    \end{tabular}
    \caption{Polynomial interpolation of the integral of $u(x)$ along a cell, as used in the tranverse flux reconstruction method. In this method, the polynomial $Q_i(x)$ is reconstructed from a neighborhood of $u(x)$ evaluated at cell midpoints. Stencils are given for polynomial degree $p=0,\ldots,7$.}
    \label{tab:transverse_integral}
\end{table}

\subsubsection{\textit{A priori} slope limiting in two dimensions}

\begin{figure}
    \centering
    \begin{subfigure}{.45\linewidth}
        \centering
        \includegraphics[width=\linewidth]{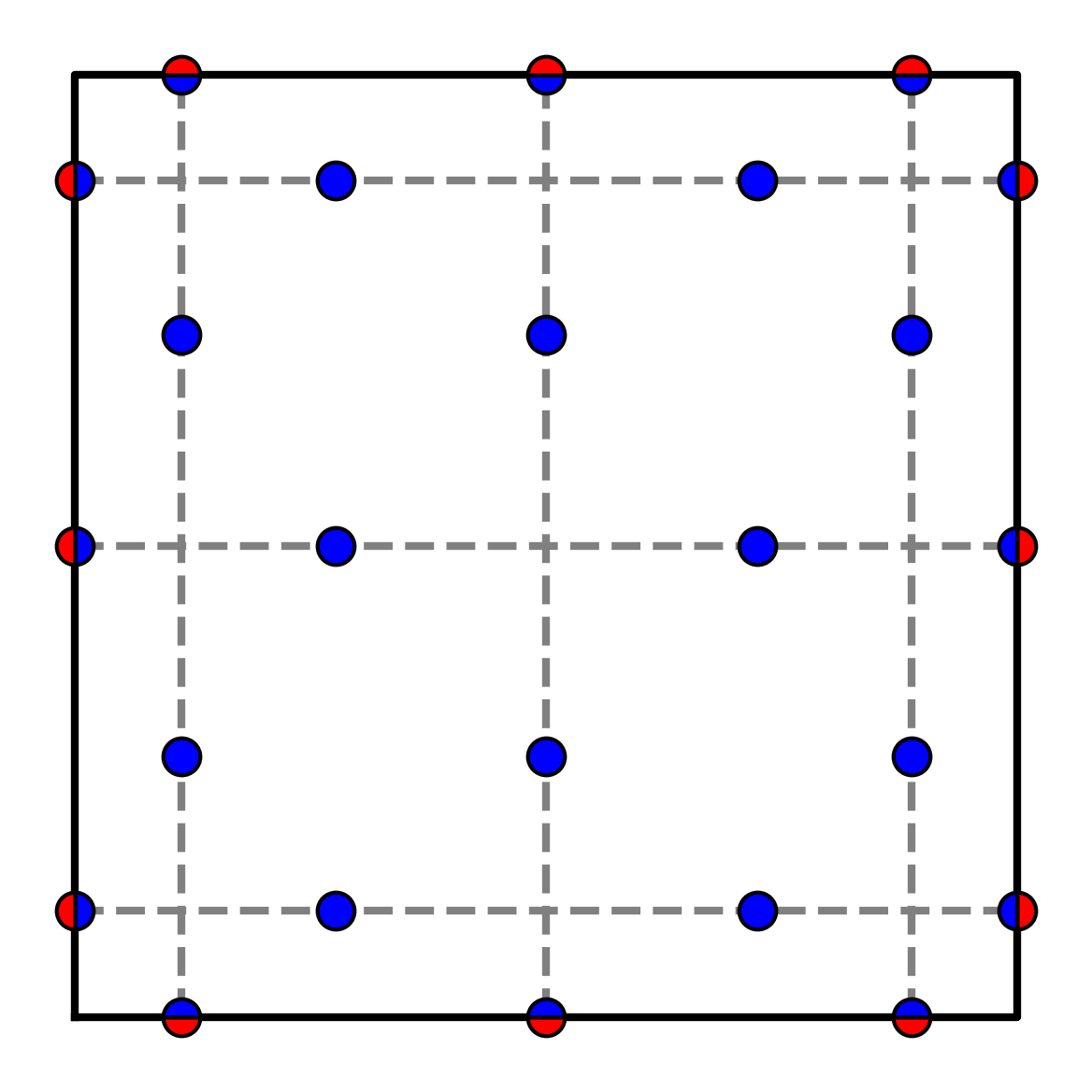}
        \caption{Gauss-Legendre set for computing fluxes (in red) and Gauss-Lobatto set for computing the limiter $\theta$ (in blue) used in the original aPrioriMPP scheme of \citep{zhang2010maximum}.}\label{fig:gaussgausslobato}
    \end{subfigure}
        \hfill
    \begin{subfigure}{.45\linewidth}
        \centering
        \includegraphics[width=\linewidth]{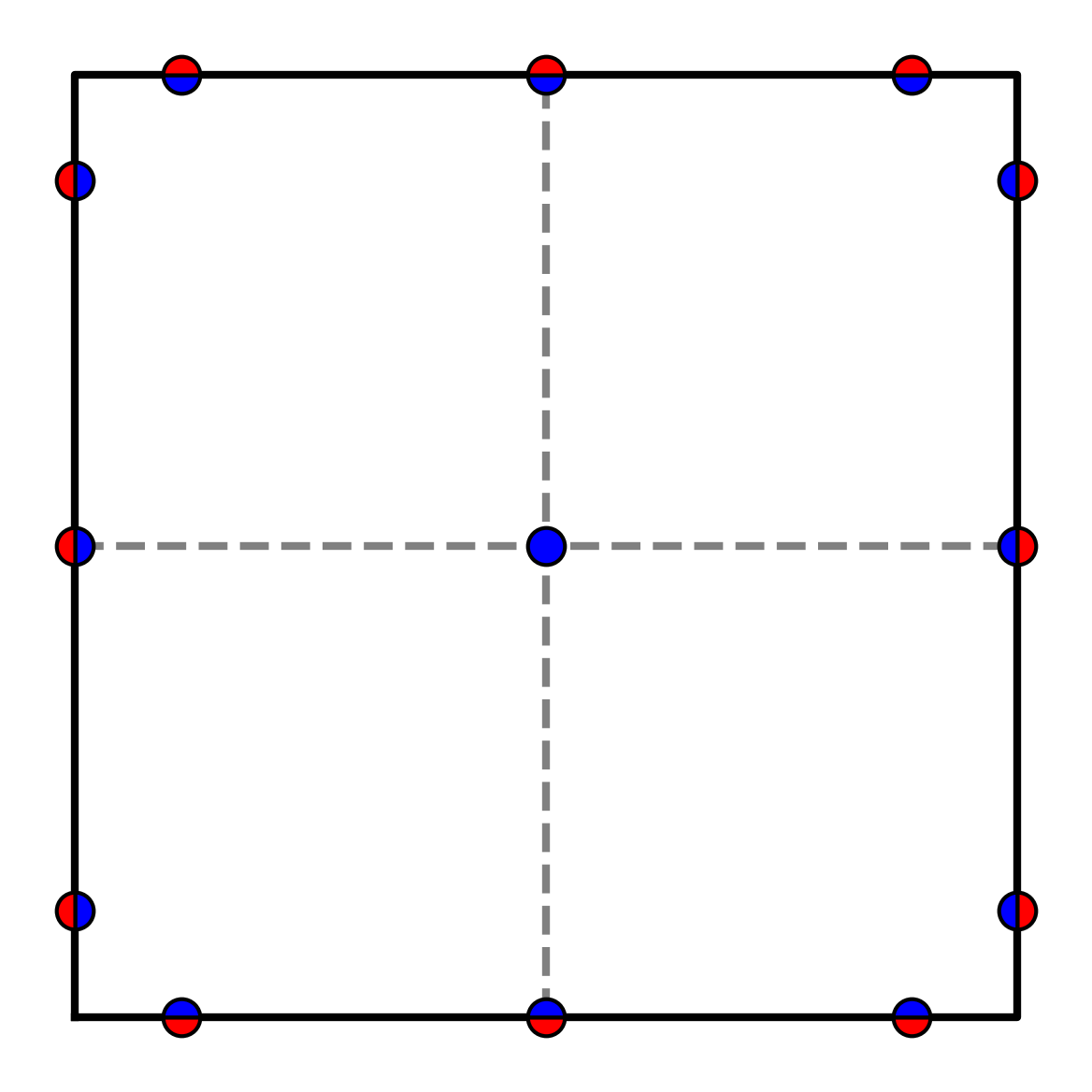}
        \caption{Gauss-Legendre set for computing fluxes (in red) and centroid set for computing the limiter $\theta$ (in blue), as used in our modified aPrioriMPP scheme.}\label{fig:gausscentroid}
    \end{subfigure}
    
    \bigskip
    
    \begin{subfigure}{.45\linewidth}
      \centering
      \includegraphics[width=\linewidth]{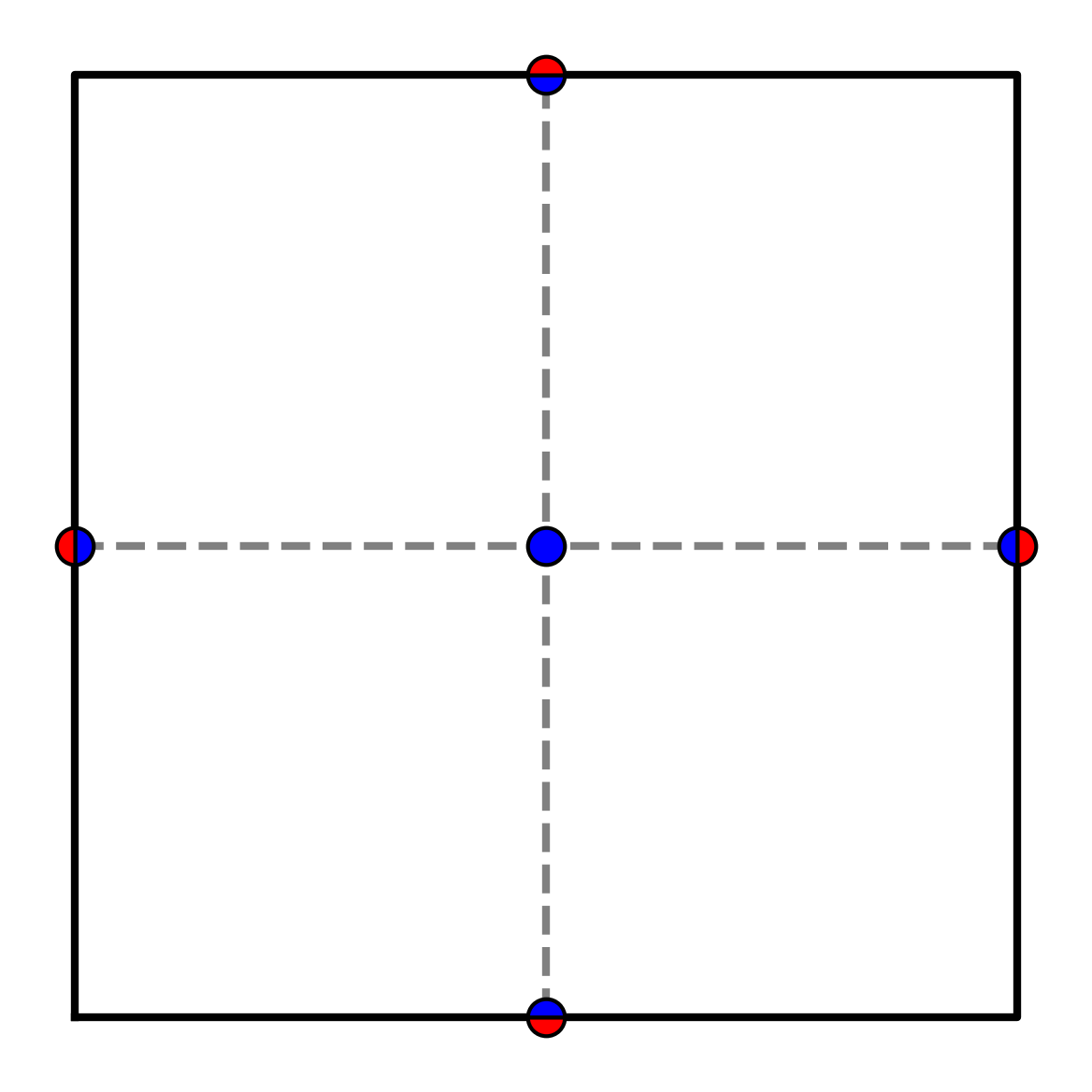}
      \caption{Transverse set for computing fluxes (in red) and centroid set for computing $\theta$ (in blue), as used in the aPrioriT scheme.}\label{fig:transversecentroid}
    \end{subfigure} 
    \caption{Cell nodal values used to compute cell fluxes (red) and the \textit{a priori} slope limiter $\theta$ (blue) for polynomial degree $p=5$.}
    \label{fig:quadratures}
\end{figure}

Zhang \& Shu's \textit{a priori} slope limited polynomial reconstruct (\cite{zhang2011maximum, zhang2010maximum}) in two-dimensions writes

\begin{equation}
    \Tilde{P}_{i,j}(x) = \theta_{i,j} \left( P_{i,j}(x) - \overline{\overline{u}}_{i,j} \right) + \overline{\overline{u}}_{i,j}, \quad \theta_{i,j} = \min \left( \left| \frac{M_{i,j} - \overline{\overline{u}}_{i,j}}{M'_{i,j} - \overline{\overline{u}}_{i,j}} \right|, \left| \frac{m_{i,j} - \overline{\overline{u}}_{i,j}}{m'_{i,j} - \overline{\overline{u}}_{i,j}} \right|, 1 \right).
\end{equation} Here, $M_{i,j}$ and $m_{i,j}$ represent the maximum and minimum, respectively, taken over the neighboring cells $\{ \overline{\overline{u}}_{i,j}, \overline{\overline{u}}_{i\pm1,j}, \overline{\overline{u}}_{i,j\pm1} \}$. As in one dimension, $M'_{i,j}$ and $m'_{i,j}$ are the maximum and minimum of the nodal interpolations of $u(x,y)$ corresponding to the $L$ Gauss-Lobatto quadrature points, but there is added complexity in two-dimensions. When the $K$-point Gauss-Legendre quadrature is used to evaluate the integral of the cell face fluxes, the $L$-points corresponding to the Gauss-Lobatto quadrature are reconstructed along each of the $K$ traces, in both the $x$- and $y$- directions (see, for example, Figure \ref{fig:gaussgausslobato}). 

To avoid the significant cost associated to this complex procedure, we introduce again the centroid set but in two dimensions. Here, we compute $M'_{i,j}$ and $m'_{i,j}$ using the nodes along cell faces used for flux computations as well as the cell centroid, without resorting to the Gauss-Lobatto quadrature points in the cell interior. Examples of this approach are illustrated in Figures \ref{fig:gausscentroid} and \ref{fig:transversecentroid} for Gauss-Legendre and transverse flux reconstructions, respectively.

\subsubsection{\textit{A posteriori} slope limiting in two dimensions}

The second-order, slope-limited MUSCL interpolations in two dimensions are written as

\begin{equation}
    \begin{split}
        \Tilde{u}_{i-\onehalf,j}^+ = \overline{\overline{u}}_{i,j} - \frac{1}{2} \Tilde{S}_{i,j}^x, \quad
        \Tilde{u}_{i+\onehalf,j}^- = \overline{\overline{u}}_{i,j} + \frac{1}{2} \Tilde{S}_{i,j}^x \\
        \Tilde{u}_{i,j-\onehalf}^+ = \overline{\overline{u}}_{i,j} - \frac{1}{2} \Tilde{S}_{i,j}^y, \quad
        \Tilde{u}_{i,j+\onehalf}^- = \overline{\overline{u}}_{i,j} + \frac{1}{2} \Tilde{S}_{i,j}^y, \\
    \end{split}
\end{equation} where the limited slopes are defined in both the $x$ and $y$ directions:

\begin{equation}
    \begin{split}
        \Tilde{S}_{i,j}^x=\text{lim}(\overline{\overline{u}}_{i-1,j}, \overline{\overline{u}}_{i,j}, \overline{\overline{u}}_{i+1,j}) \frac{\overline{\overline{u}}_{i+1,j}-\overline{\overline{u}}_{i-1,j}}{2} \\
        \Tilde{S}_{i,j}^y=\text{lim}(\overline{\overline{u}}_{i,j-1}, \overline{\overline{u}}_{i,j}, \overline{\overline{u}}_{i,j+1}) \frac{\overline{\overline{u}}_{i,j+1}-\overline{\overline{u}}_{i,j-1}}{2}. \\
    \end{split}
    \label{general_limited_slopes_xy}
\end{equation}

The minmod and moncen slope limiters act independently in both directions. However, neither slope limiter used in this way results in a strictly MPP MUSCL scheme. Therefore, we introduce a third slope limiter: the positivity-preserving (PP) slope limiter, which was first presented by \cite{suresh2000positivity} as a means to allow MUSCl-Hancock to be positivity-preserving in two-dimensions.

The PP slope limiter neighbor differences are defined as

\begin{equation}
    V_{min}=\min \begin{bmatrix}
    \overline{\overline{u}}_{i-1,j+1} - \overline{\overline{u}}_{i,j} & \overline{\overline{u}}_{i,j+1} - \overline{\overline{u}}_{i,j} & \overline{\overline{u}}_{i+1,j+1} - \overline{\overline{u}}_{i,j} \\
    \overline{\overline{u}}_{i-1,j} - \overline{\overline{u}}_{i,j} & -\epsilon_{pp} & \overline{\overline{u}}_{i+1,j} - \overline{\overline{u}}_{i,j} \\
    \overline{\overline{u}}_{i-1,j-1} - \overline{\overline{u}}_{i,j} & \overline{\overline{u}}_{i,j-1} - \overline{\overline{u}}_{i,j} & \overline{\overline{u}}_{i+1,j-1} - \overline{\overline{u}}_{i,j} \\
    \end{bmatrix}
\end{equation} and

\begin{equation}
    V_{max}=\max \begin{bmatrix}
    \overline{\overline{u}}_{i-1,j+1} - \overline{\overline{u}}_{i,j} & \overline{\overline{u}}_{i,j+1} - \overline{\overline{u}}_{i,j} & \overline{\overline{u}}_{i+1,j+1} - \overline{\overline{u}}_{i,j} \\
    \overline{\overline{u}}_{i-1,j} - \overline{\overline{u}}_{i,j} & \epsilon_{pp} & \overline{\overline{u}}_{i+1,j} - \overline{\overline{u}}_{i,j} \\
    \overline{\overline{u}}_{i-1,j-1} - \overline{\overline{u}}_{i,j} & \overline{\overline{u}}_{i,j-1} - \overline{\overline{u}}_{i,j} & \overline{\overline{u}}_{i+1,j-1} - \overline{\overline{u}}_{i,j} \\
    \end{bmatrix},
\end{equation} where $\epsilon_{pp}=1\times10^{-20}$. The limited slopes (\ref{general_limited_slopes_xy}) become

\begin{equation}
    \begin{split}
        \Tilde{S}_{i,j}^x = \min(1,V) \frac{\overline{\overline{u}}_{i+1,j} - \overline{\overline{u}}_{i-1,j}}{2} \\
        \Tilde{S}_{i,j}^y = \min(1,V) \frac{\overline{\overline{u}}_{i,j+1} - \overline{\overline{u}}_{i,j-1}}{2}
    \end{split}
\end{equation} where

\begin{equation}
    V=\frac{4 \cdot \min(|V_{min}|,|V_{max}|)}{|\overline{\overline{u}}_{i+1,j} - \overline{\overline{u}}_{i-1,j}| + |\overline{\overline{u}}_{i,j+1} - \overline{\overline{u}}_{i,j-1}|}.
\end{equation}

With the fallback scheme established in two dimensions, we define the convex blending of corrected fluxes presented by Vilar \& Abgrall \cite{vilar2024posteriori}. Let the high-order face flux integrals be written

\begin{equation}
    \hat{F}_{i \pm \onehalf,j}^p = \int_{y_{j-\onehalf}}^{y_{j+\onehalf}} \hat{f}[ u_{i \pm \onehalf,j}^-(\eta) , u_{i \pm \onehalf,j}^+(\eta) ] d\eta
\end{equation} and

\begin{equation}
    \hat{G}_{i, j \pm \onehalf}^p = \int_{x_{i-\onehalf}}^{x_{i+\onehalf}} \hat{g}[ u_{i,j \pm \onehalf}^-(\xi) , u_{i,j \pm \onehalf}^+(\xi) ] d\xi.
\end{equation} Similarly, let the fallback face flux integrals be written

\begin{equation}
    \hat{F}_{i \pm \onehalf, j}^1 = \hat{f}(\Tilde{u}_{i \pm \onehalf, j}^-, \Tilde{u}_{i \pm \onehalf, j}^+)
\end{equation} and 

\begin{equation}
    \hat{G}_{i, j \pm \onehalf}^1 = \hat{g}(\Tilde{u}_{i, j \pm \onehalf}^-, \Tilde{u}_{i, j \pm \onehalf}^+).
\end{equation} An integral quadrature is not needed because the cell face midpoint is already a second-order approximation (see Table \ref{tab:transverse_integral}). Then, the blending formula (\ref{convexblending}) becomes 

\begin{equation}
    \begin{split}
        \hat{F}_{i - \onehalf,j}^p \leftarrow \max(\phi_{i-1,j}, \phi_{i,j}) (\hat{F}^1_{i-\onehalf, j} - \hat{F}_{i-\onehalf, j}^p) + \hat{F}_{i-\onehalf, j}^p \\
        \hat{F}_{i + \onehalf,j}^p \leftarrow \max(\phi_{i,j}, \phi_{i+1,j}) (\hat{F}^1_{i+\onehalf, j} - \hat{F}_{i+\onehalf, j}^p) + \hat{F}_{i+\onehalf, j}^p \\
        \hat{G}_{i,j - \onehalf}^p \leftarrow \max(\phi_{i,j-1}, \phi_{i,j}) (\hat{G}^1_{i,j-\onehalf} - \hat{G}_{i,j-\onehalf}^p) + \hat{G}_{i, j-\onehalf}^p \\
        \hat{G}_{i, j + \onehalf}^p \leftarrow \max(\phi_{i,j}, \phi_{i,j+1}) (\hat{G}^1_{i,j+\onehalf} - \hat{G}_{i,j+\onehalf}^p) + \hat{G}_{i,j-\onehalf}^p
    \end{split}
\end{equation} in two dimensions.

\begin{figure}
    \centering
    \begin{subfigure}{.45\linewidth}
        \centering
        \includegraphics[width=\linewidth]{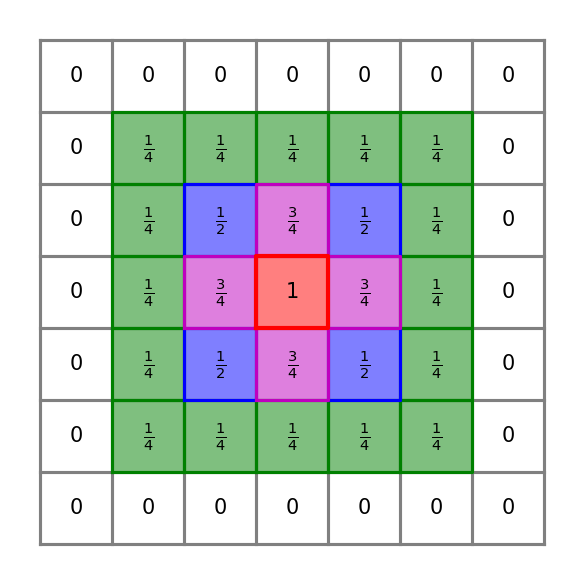}
        \caption{One troubled cell.}\label{fig:image1}
    \end{subfigure}
        \hfill
    \begin{subfigure}{.45\linewidth}
        \centering
        \includegraphics[width=\linewidth]{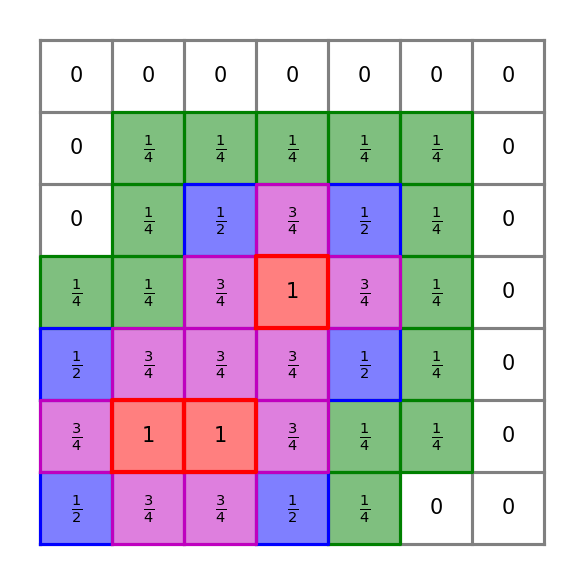}
        \caption{Three troubled cells.}\label{fig:image12}
    \end{subfigure}
    \caption{The value of $\phi$ used for the convex blending of corrected fluxes, shown across a region comprising 49 cells with troubled cells highlighted in red.}
    \label{fig:blending2d}
\end{figure}

The blending parameter $\phi$ is computed at each cell via the following procedure introduced by Vilar \& Abgrall: If cell $i,j$ is troubled, set $\phi_{i,j}=1$. If cell $i,j$ is not troubled and it shares a face with a troubled cell, set $\phi_{i,j}=\frac{3}{4}$. If cell $i,j$ is not troubled and it shares a corner with a troubled cell, but not a face, set $\phi_{i,j}=\frac{1}{2}$. If cell $i,j$ is not troubled and its Euclidean distance from the nearest troubled cell is 2 cell units, set $\phi_{i,j} = \frac{1}{4}$. If cell $i,j$ is not troubled and its Euclidean distance from the nearest troubled cell is greater than 2 cell units, set $\phi_{i,j}=0$. See Figure \ref{fig:blending2d} for examples.

\subsubsection{Smooth extrema detection in two dimensions}

The one-dimensional smoothness indicator $\alpha$, as defined in (\ref{smooth_extrema_detection}), is computed at each cell in both the $x$- and $y$-directions. Let $\alpha^x$ and $\alpha^y$ denote the indicators computed in these directions. Then, cell $i,j$ is considered a smooth extremum only if

\begin{equation}
    \min(\alpha^x_{i-1,j}, \alpha^x_{i,j}, \alpha^x_{i+1,j}) = \min(\alpha^y_{i,j-1}, \alpha^y_{i,j}, \alpha^y_{i,j+1}) = 1.
\end{equation}

\section{Summary of numerical schemes}
\label{sec:numerical_schemes}

We have seen a variety of options for slope limiters as modifications to the high-order finite volume semi-discrete scheme. We summarize in this section the semi-discrete schemes that are studied via numerical experiments in Section \ref{sec:numerical_results}. In all presented semi-discrete schemes, the degree of the conservative interpolation polynomial $p$ is left unspecified, allowing for variation in the numerical experiments. It is important to note that the presented schemes are only semi-discrete; the fully-discrete scheme is only defined once a Runge-Kutta integration method is chosen.

For the one-dimensional advection equation, we examine three semi-discrete schemes: aPrioriMPP, aPosteriori, and aPosterioriB. These schemes are named based on their adoption of either the \textit{a priori} or \textit{a posteriori} slope limiting paradigm, as summarized in Table \ref{tab:schemes_1d}.

In aPrioriMPP, the slope limiter $\theta$, is computed using the centroid set instead of the Gauss-Lobatto set of Zhang \& Shu \cite{zhang2011maximum, zhang2010maximum}, a substitution we find sufficient for maintaining the MPP property. We also replace Zhang \& Shu's reduced CFL factor $C_{\text{MPP}}$ \cite{zhang2011maximum, zhang2010maximum} with the adaptive time-step size of Huang \textit{et al} \cite{huang2023high}. The initial time-step size is deduced from an initial CFL factor of 0.8. This way, aPrioriMPP achieves more competitive speed performance.

The \textit{a posteriori} slope limiting semi-discrete schemes aPosteriori and aPosterioriB use the moncen slope limiter in their MUSCL fallback scheme. aPosterioriB features the convex blending of revised fluxes introduced by Vilar \& Abgrall \cite{vilar2024posteriori}.

In the context of the two-dimensional advection equation, detailed semi-discrete schemes are outlined in Table \ref{tab:schemes_2d}. An essential specification for these schemes is the flux integral method. Specifically, aPrioriMPP constructs cell face fluxes with Gauss-Legendre quadrature, while aPrioriT uses a transverse flux reconstruction. We disable the adaptive time-step size for aPrioriT due to maximum principle violations observed with this combination (even when the reduced CFL factor $C_{\text{MPP}}$ is used).

In two dimensions, we observe that \textit{a posteriori} schemes exhibit similar maximum principle violations regardless of flux reconstruction method, so we opt for the cheaper option (transverse) in aPosteriori and aPosterioriB. PP, the two-dimensional slope limiter of Suresh \cite{suresh2000positivity} is chosen as the slope limiter of the MUSCL fallback scheme in these cases.

Smooth extrema detection is enabled across all semi-discrete schemes in both one and two dimensions.

\begin{table}[]
    \centering
    \begin{tabular}{c|ccccc}
        \toprule
         Scheme & $\theta_i$ set & $C$ & Adaptive $\Delta t$ & Fallback limiter & Blending \\ 
         \midrule
         aPrioriMPP & Centroid & 0.8 & Yes & --- & --- \\
         \cline{1-6}
         aPosteriori & --- & 0.8 & No & moncen & No \\
         \cline{1-6}
         aPosterioriB & --- & 0.8 & No & moncen & Yes \\
         \bottomrule
    \end{tabular}
    \caption{Summary of three slope-limited, semi-discrete schemes for solving the one-dimensional advection equation.}
    \label{tab:schemes_1d}
\end{table}

\begin{table}[]
    \centering
    \begin{tabular}{c|ccccccc}
        \toprule
         Scheme & Flux reconstruction & $\theta_{ij}$ set & $C$ & Adaptive $\Delta t$ & Fallback limiter & Blending \\ 
         \midrule
         aPrioriMPP & Gauss-Legendre & Centroid & 0.8 & Yes & --- & --- \\
         \cline{1-7}
         aPrioriT & Transverse & Centroid & 0.8 & No & --- & --- \\
         \cline{1-7}
         aPosteriori & Transverse & --- & 0.8 & No & PP & No \\
         \cline{1-7}
         aPosterioriB & Transverse & --- & 0.8 & No & PP & Yes \\
         \bottomrule
    \end{tabular}
    \caption{Summary of four slope-limited, semi-discrete schemes for solving the two-dimensional advection equation.}
    \label{tab:schemes_2d}
\end{table}

\section{Numerical results}
\label{sec:numerical_results}

In this section, we perform a series of numerical tests of our implementation of the finite volume method and the performance of the various slope limited schemes. In these tests, the linear advection equation is solved in 1D and 2D. The polynomial degree of the spatial interpolation polynomial is given by $p$, the number of cells in one direction is represented by $N$, and for 2D solutions, an array of $N \times N$ cells is used. $h=\frac{L}{N}$ is the uniform grid spacing, where $L$ is the length of the computational domain.

The maximum principle of a linear advection problem with initial condition $u(0, \mathbf{x})=u_0(\mathbf{x})$ is given by $M=\max_\mathbf{x} u_0(\mathbf{x})$ and $m=\min_\mathbf{x} u_0(\mathbf{x})$. To indicate the presence of a maximum principle violation in our numerical solutions, we define

\begin{equation}
    \delta=\min (\delta^-, \delta^+)
    \label{delta}
\end{equation} and

\begin{align*}
    \delta^- = \min_n \min_{ij} (\overline{\overline{u}}_{ij}^n - m), \quad \delta^+ = \min_n \min_{ij} (M - \overline{\overline{u}}_{ij}^n),
\end{align*} where the indices $i$ and $j$ cover the entire computational domain and the index $n$ covers all time-steps (see \cite{kuzmin2022bound}). Negative values of $\delta$ imply a violation of the maximum principle. We consider numerical solutions with $\delta >$ -1E-10 to sufficiently preserve the maximum principle and refer to such solutions as \textit{approximately MPP}.

\subsection{One-dimensional advection of the composite profile}

We perform a 1D advection test using the classic composite profile for $\bar{u}_0(x)$  \cite{kuzmin2022bound, krivodonova2007limiters, jiang1996efficient}, with a velocity $a=1$, and a periodic region $x\in[0,1]$. The problem is solved with aPrioriMPP, aPosteriori, and aPosterioriB schemes up to $t=1$ for various polynomial degrees $p$ and Runge-Kutta methods.

Results are shown in Table \ref{tab:composite}. aPrioriMPP produces very good MPP numerical solutions for all given values of $p$ and all SSP Runge-Kutta methods, as well as RK4. In contrast, \textit{a posteriori} solutions exhibit violations in all cases. Using RK4 instead of SSPRK3 results in smaller violations for the \textit{a posteriori} slope limited schemes. These violations are made smaller still through the use of convex blending for the revised fluxes, in some cases by several orders of magnitude. We observe that the magnitude of the maximum principle violations of the \textit{a posteriori} slope limited schemes do not strictly decrease with $N$ or the number of time-steps. In general, they do decrease with increasing $p$.

\begin{table}[]
    \centering
    \begin{tabular}{llrrr}
    \toprule
    & & \multicolumn{3}{c}{$\delta$} \\
    \cline{3-5}
     $p$ & Integrator & aPrioriMPP & aPosteriori & aPosterioriB \\
    \midrule
    1 & SSPRK2 & \textbf{-2.22E-18} & -1.04E-02 & -4.57E-03 \\
    \cline{1-5}
    2 & SSPRK3 & \textbf{-2.82E-19} & -6.52E-03 & -8.70E-04 \\
    \cline{1-5}
    \multirow[t]{2}{*}{3} & SSPRK3 & \textbf{-2.12E-11} & -7.85E-03 & -2.05E-04 \\
     & RK4 & \textbf{-1.10E-11} & -6.44E-05 & -5.37E-05 \\
    \cline{1-5}
    \multirow[t]{2}{*}{4} & SSPRK3 & \textbf{-1.58E-12} & -6.83E-03 & -2.83E-04 \\
     & RK4 & \textbf{-3.91E-11} & -4.13E-05 & -6.43E-05 \\
    \cline{1-5}
    \multirow[t]{2}{*}{5} & SSPRK3 & \textbf{-5.43E-11} & -7.91E-03 & -2.08E-04 \\
     & RK4 & \textbf{-3.27E-12} & -1.00E-04 & -1.40E-08 \\
    \cline{1-5}
    \multirow[t]{2}{*}{6} & SSPRK3 & \textbf{-1.43E-11} & -7.46E-03 & -1.97E-04 \\
     & RK4 & \textbf{-1.50E-13} & -5.86E-07 & -2.38E-07 \\
    \cline{1-5}
    \multirow[t]{2}{*}{7} & SSPRK3 & \textbf{-4.33E-14} & -7.70E-03 & -2.64E-04 \\
     & RK4 & \textbf{-5.70E-11} & -3.35E-04 & -1.81E-06 \\
    \cline{1-5}
    \bottomrule
    \end{tabular}
    \caption{Maximum principle violation $\delta$ of three types of schemes after solving the 1D advection of the composite profile up to $t=1$. The test is repeated for various Runge-Kutta methods and polynomial degree $p$. Violations smaller in magnitude than -1E-10 are typed in boldface font.}
    \label{tab:composite}
\end{table}

\begin{figure}
    \centering
    \includegraphics[width=\textwidth]{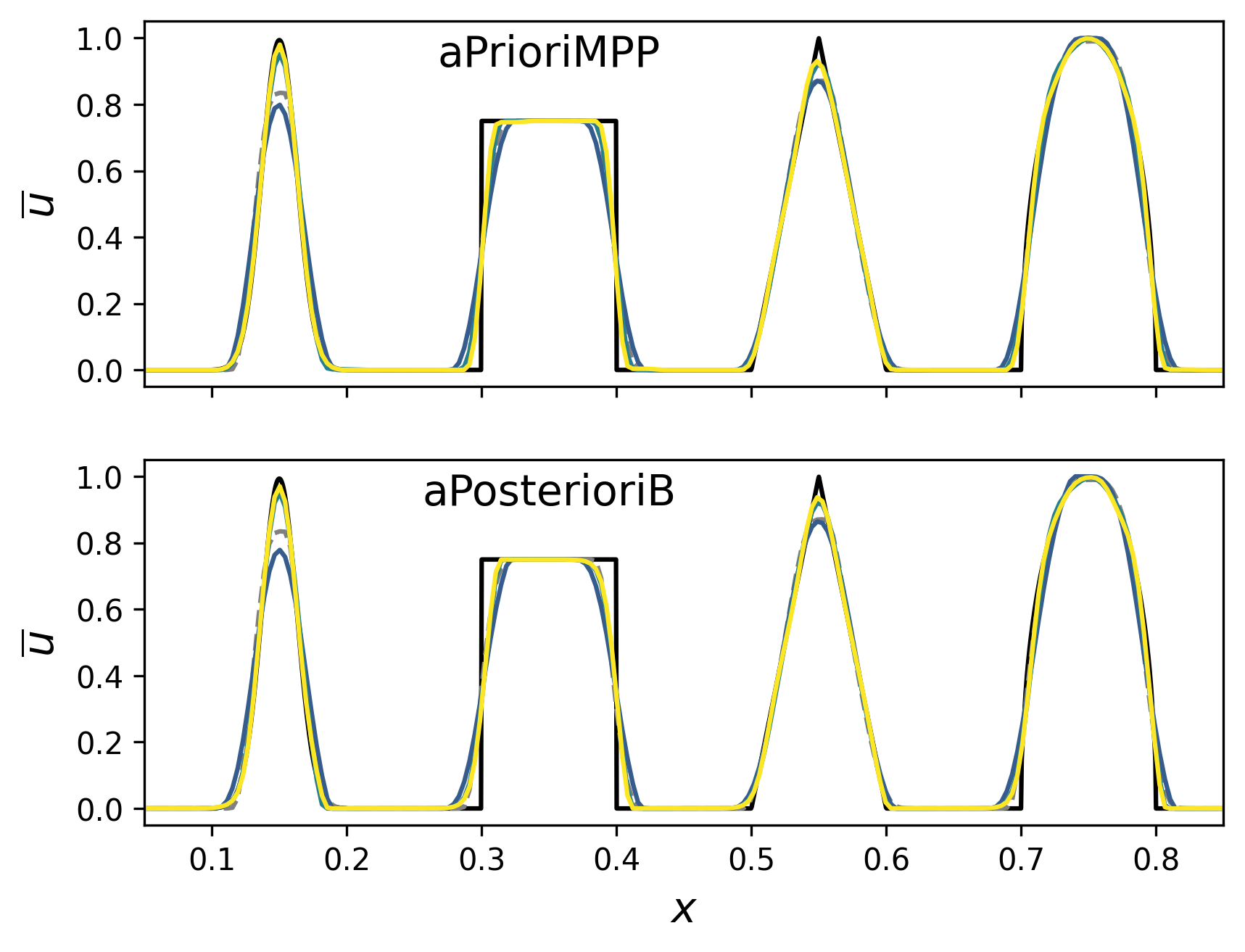}
    \caption{Snapshots of the numerical solution to the advection of the composite profile at $t=1$. Results are shown for the aPrioriMPP and aPosterioriB schemes and for polynomial degrees $p=2$ (dark blue), $p=3$ (light blue), and $p=7$ (yellow). The aPrioriMPP schemes use SSPRK3 for all results shown while the aPosterioriB schemes use SSPRK3 for $p=2$ and RK4 for $p>2$. The numerical solution of second-order MUSCL-Hancock is shown in dashed grey for reference. All results shown have a resolution of $N=256$ cells. Maximum principle violations are not observed in the aPrioriMPP and MUSCL-Hancock solutions while they are observed for the aPosterioriB solutions.}
    \label{fig:composite1}
\end{figure}

\begin{figure}
    \centering
    \includegraphics[width=\textwidth]{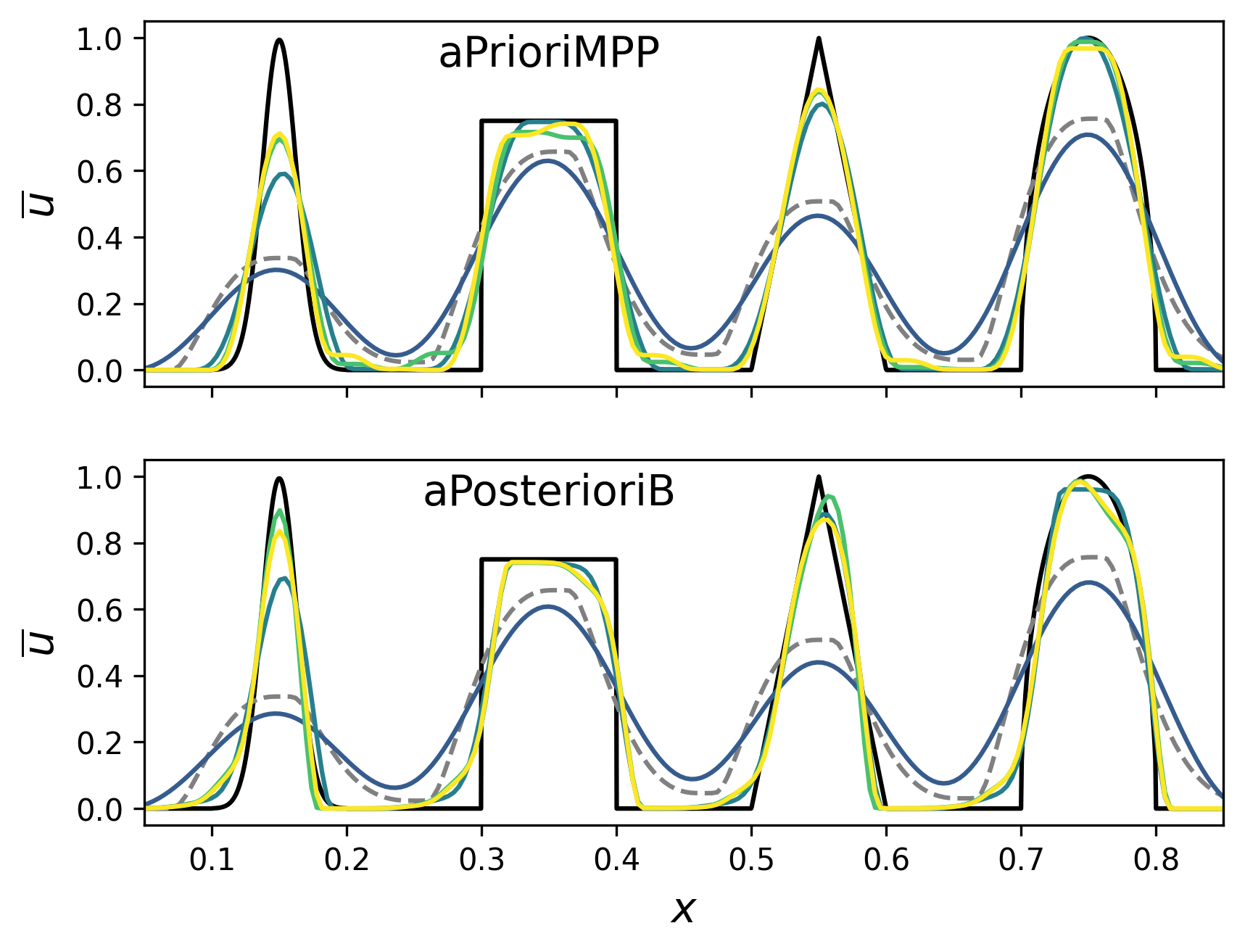}
    \caption{Snapshots of the numerical solution to the advection of the composite profile at $t=100$. Results are shown for the aPrioriMPP and aPosterioriB schemes and for polynomial degrees $p=2$ (dark blue), $p=3$ (light blue), $p=5$ (green), and $p=7$ (yellow). The aPrioriMPP schemes use SSPRK3 for all results shown while the aPosterioriB schemes use SSPRK3 for $p=2$ and RK4 for $p>2$. The numerical solution of second-order MUSCL-Hancock is shown in dashed grey for reference. All results shown have a resolution of $N=256$ cells. Maximum principle violations are not observed in the aPrioriMPP and MUSCL-Hancock solutions while they are observed for the aPosterioriB solutions.}
    \label{fig:composite100}
\end{figure}

Figures \ref{fig:composite1} and \ref{fig:composite100} shows our numerical solutions at different times generated by aPrioriMPP and aPosterioriB schemes as well as second-order MUSCL-Hancock. For $p>2$, aPrioriMPP is solved with SSPRK3 since RK4 is found to result in excessive numerical diffusion. Meanwhile, aPosterioriB uses RK4 for $p>3$ and RK4 since this combination exhibits smaller violations of the maximum principle, as seen in Table \ref{tab:composite}.

After one period of advection, the numerical solutions of aPrioriMPP and aPosterioriB appear quite similar when $p$ is the same. Those with $p>2$ are of particularly good quality at this time, but the key difference remains that aPrioriMPP maintains the maximum principle for all $p$ while aPosterioriB results in significant violations.

After 100 periods, the higher-$p$ ($p>2$) schemes show greater resilience to numerical diffusion than their lower-order counterparts, as they better preserve the initial composite profile. Numerical artifacts become apparent as $p$ increases at this longer time, with aPrioriMPP showing slightly more pronounced artifacts. MUSCL-Hancock outperforms our third-order schemes, but shows greater numerical diffusion than our higher-order implementations.

Comparing our results for aPrioriMPP at $p=3$ and $p=5$ to a similar test conducted by Kuzmin \textit{et al} with their own MPP schemes \cite{kuzmin2022bound}, we observe significantly less numerical diffusion and artifacts in our implementation. This is because Kuzmin \textit{et al} use a fourth-order, five-stage and sixth-order, seven-stage Runge-Kutta method for $p=3$ and $p=5$, respectively \cite{kuzmin2022bound}, while we use SSPRK3 with only three stages in these cases. Granted, this makes their fully-discrete schemes truly high-order, while ours are capped at third-order by SSPRK3; The benefits of very-high-order finite volume schemes are highlighted in a later experiment. However, there is a clear disadvantage of the very-high-order time integration for problems with non-smooth initial data, such as the composite profile.

As is shown in later experiments, the \textit{a priori} limiting schemes of both Zhang \& Shu \cite{zhang2011maximum, zhang2010maximum} and Kuzmin \textit{et al} \cite{kuzmin2022bound} are much more conservative than the \textit{a posteriori} limiting schemes in the magnitude of slopes they permit in the presence of discontinuities. Consequently, the \textit{a priori} limiting schemes have a tendency to produce numerical artifacts in these regions. These artifacts are exacerbated by the large number of steps needed for long time integration, so the high-order time integrators actually worsen the quality of discontinuous solutions. Thus, we find very different performance for the \textit{a priori} limiting schemes between smooth and discontinuous solutions. 

\subsection{Two-dimensional advection of a sine wave}

We conduct the 2D advection of a smooth solution to check the order of accuracy of our implementation of the finite volume method. The setup

\begin{equation}
    \overline{\overline{u}}_0(x,y) = \sin (2 \pi (x+y)), \quad v_x=2, \quad v_y=1,
\end{equation} is evolved in a periodic box $x, y \in [0,1]$. The $L_1$ norm of the error of the numerical solution is computed at $t=1$ (one period of advection) using the formula

\begin{equation}
    E_1= h^2 \sum_{i,j} |\overline{\overline{u}}_{ij} - \overline{\overline{u}}_{0, ij}|,
\end{equation} where the indices $i$ and $j$ cover the $x$- and $y$-components of our computational domain.

We vary $N$ and $p$, ensuring the order of accuracy of the chosen Runge-Kutta method is no less than the order of accuracy of the spatial discretization. For $p>5$, where we lack a higher-than-sixth order Runge-Kutta method, we adopt Vilar's \cite{vilar2024posteriori} approach of reducing $\Delta t$ to synthesize higher-order solutions. With this approach, $C$ is determined by the formula:

\begin{equation}
    C_p =
    \begin{cases}
    0.8 \left( \frac{h}{L} \right)^{\frac{p-q}{q+1}} & \text{if } p>q\\
    0.8             & \text{otherwise},
    \end{cases}
\label{time_step_modification}
\end{equation} where $q$ is the polynomial degree of the temporal discretization.

The smooth sine wave problem is solved using aPrioriMPP, which uses a Gauss-Legendre flux quadrature, and aPrioriT, which uses a transverse flux reconstruction. The adaptive time-step size of aPrioriMPP is disabled since the CFL factor $C$ is chosen according to (\ref{time_step_modification}). We also skip the check for overshoots and undershoots of the maximum principle  in the smooth extrema detection routine, since preserving a strict maximum principle is not the goal of this experiment.

\begin{figure}
    \centering
    \includegraphics[width=\textwidth]{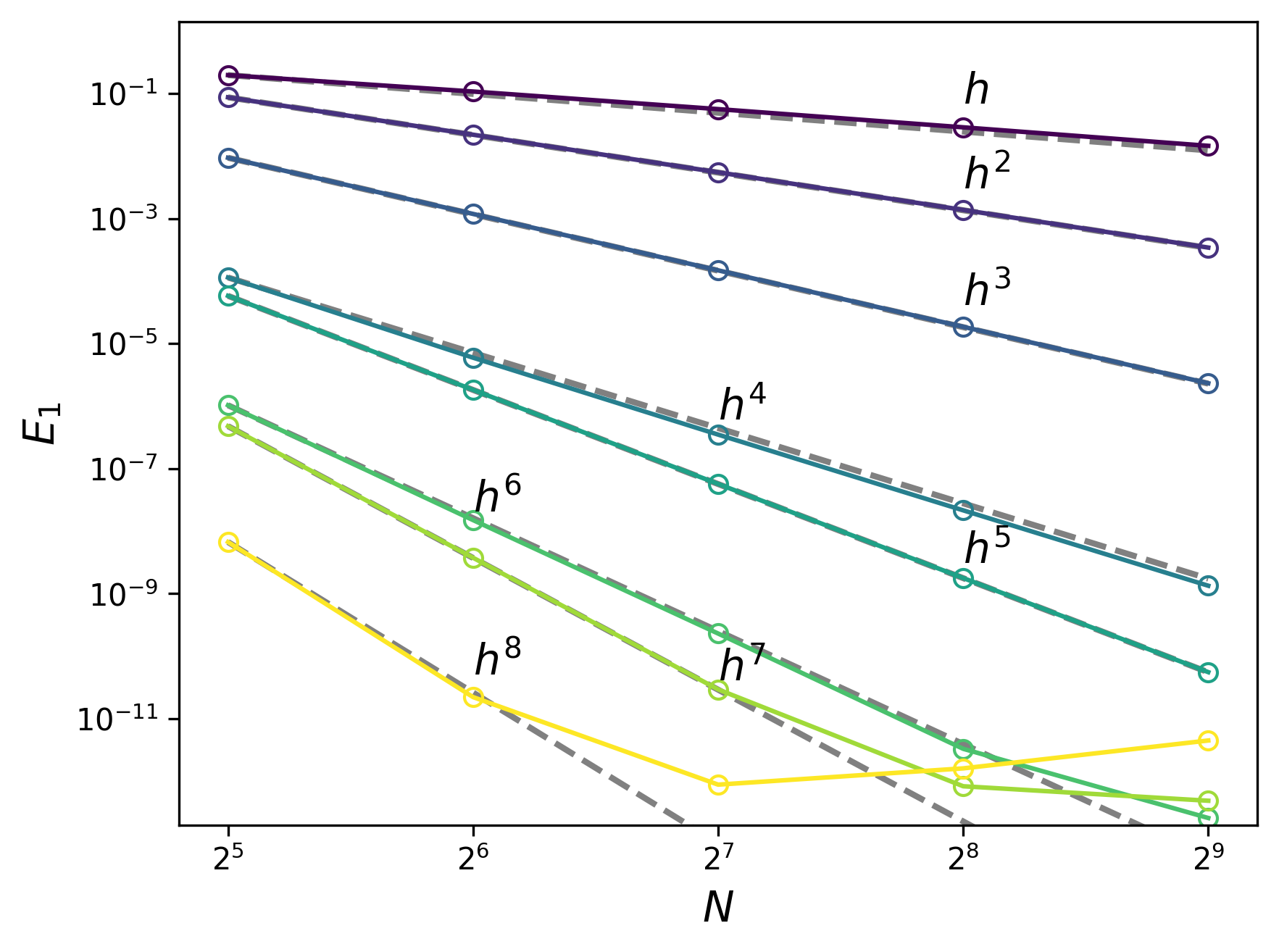}
    \caption{The $L_1$ error of aPrioriMPP schemes at $t=1$ for the 2D advection of a sine wave with varying resolution $N$. Results are shown for polynomial degree $p$ varying from 0 to 7, shaded linearly from blue to yellow. For this experiment, the adaptive time-step size is disabled and $C$ is chosen according to (\ref{time_step_modification}).}
    \label{fig:convergence}
\end{figure}

The results are summarized in Table \ref{tab:convergence} with the error convergence of the aPrioriMPP schemes depicted in Figure \ref{fig:convergence}. It is shown that the errors of the numerical solutions converge as $N$ increases at the expected rates, consistent with the designed order of accuracy for each scheme, until a precision floor is reached at around 1E-12. This is true for both the Gauss-Legendre quadrature and transverse flux reconstruction. Both flux reconstructions produce the same error when $p$ is even or less than 3. However, for $p=3,5,7$,  the transverse reconstruction yields a smaller error, sometimes up to 40\% lower. This difference is hardly visible on the log scale of the convergence study.

We emphasize the effectiveness of smooth extrema detection in this experiment. Despite the fact that slope limiting was enabled, the high-order solution was not contaminated with low-order approximations because smooth extrema detection in this case always disables the limiter.

\begin{table}[]
    \centering
    \begin{tabular}{llllllll}
	\toprule
	 &  &  &  & $E^1_\text{GL}$ & EOC & $E^1_T$ & EOC \\
	$p$ & integrator & $N$ & $C$ &  &  &  &  \\
	\midrule
	\multirow[t]{5}{*}{0} & \multirow[t]{5}{*}{Euler} & 32 & 0.80 & 1.97E-01 & --- & 1.97E-01 & --- \\
	\cline{3-8}
	 &  & 64 & 0.80 & 1.08E-01 & 0.873 & 1.08E-01 & 0.873 \\
	\cline{3-8}
	 &  & 128 & 0.80 & 5.63E-02 & 0.935 & 5.63E-02 & 0.935 \\
	\cline{3-8}
	 &  & 256 & 0.80 & 2.88E-02 & 0.967 & 2.88E-02 & 0.967 \\
	\cline{3-8}
	 &  & 512 & 0.80 & 1.46E-02 & 0.983 & 1.46E-02 & 0.983 \\
	\cline{1-8} \cline{2-8} \cline{3-8}
	\multirow[t]{5}{*}{1} & \multirow[t]{5}{*}{SSPRK2} & 32 & 0.80 & 8.69E-02 & --- & 8.69E-02 & --- \\
	\cline{3-8}
	 &  & 64 & 0.80 & 2.19E-02 & 1.986 & 2.19E-02 & 1.986 \\
	\cline{3-8}
	 &  & 128 & 0.80 & 5.49E-03 & 1.998 & 5.49E-03 & 1.998 \\
	\cline{3-8}
	 &  & 256 & 0.80 & 1.37E-03 & 1.999 & 1.37E-03 & 1.999 \\
	\cline{3-8}
	 &  & 512 & 0.80 & 3.43E-04 & 2.000 & 3.43E-04 & 2.000 \\
	\cline{1-8} \cline{2-8} \cline{3-8}
	\multirow[t]{5}{*}{2} & \multirow[t]{5}{*}{SSPRK3} & 32 & 0.80 & 9.38E-03 & --- & 9.38E-03 & --- \\
	\cline{3-8}
	 &  & 64 & 0.80 & 1.19E-03 & 2.984 & 1.19E-03 & 2.984 \\
	\cline{3-8}
	 &  & 128 & 0.80 & 1.48E-04 & 2.997 & 1.48E-04 & 2.997 \\
	\cline{3-8}
	 &  & 256 & 0.80 & 1.86E-05 & 2.999 & 1.86E-05 & 2.999 \\
	\cline{3-8}
	 &  & 512 & 0.80 & 2.32E-06 & 3.000 & 2.32E-06 & 3.000 \\
	\cline{1-8} \cline{2-8} \cline{3-8}
	\multirow[t]{5}{*}{3} & \multirow[t]{5}{*}{RK4} & 32 & 0.80 & 1.14E-04 & --- & 9.48E-05 & --- \\
	\cline{3-8}
	 &  & 64 & 0.80 & 5.95E-06 & 4.261 & 4.28E-06 & 4.469 \\
	\cline{3-8}
	 &  & 128 & 0.80 & 3.50E-07 & 4.086 & 2.34E-07 & 4.193 \\
	\cline{3-8}
	 &  & 256 & 0.80 & 2.15E-08 & 4.023 & 1.41E-08 & 4.057 \\
	\cline{3-8}
	 &  & 512 & 0.80 & 1.34E-09 & 4.006 & 8.70E-10 & 4.015 \\
	\cline{1-8} \cline{2-8} \cline{3-8}
	\multirow[t]{5}{*}{4} & \multirow[t]{5}{*}{RK6} & 32 & 0.80 & 5.81E-05 & --- & 5.81E-05 & --- \\
	\cline{3-8}
	 &  & 64 & 0.80 & 1.82E-06 & 4.994 & 1.82E-06 & 4.995 \\
	\cline{3-8}
	 &  & 128 & 0.80 & 5.70E-08 & 4.999 & 5.70E-08 & 4.999 \\
	\cline{3-8}
	 &  & 256 & 0.80 & 1.78E-09 & 5.000 & 1.78E-09 & 5.000 \\
	\cline{3-8}
	 &  & 512 & 0.80 & 5.57E-11 & 5.000 & 5.57E-11 & 5.000 \\
	\cline{1-8} \cline{2-8} \cline{3-8}
	\multirow[t]{5}{*}{5} & \multirow[t]{5}{*}{RK6} & 32 & 0.80 & 1.03E-06 & --- & 8.56E-07 & --- \\
	\cline{3-8}
	 &  & 64 & 0.80 & 1.50E-08 & 6.102 & 1.18E-08 & 6.180 \\
	\cline{3-8}
	 &  & 128 & 0.80 & 2.31E-10 & 6.028 & 1.78E-10 & 6.053 \\
	\cline{3-8}
	 &  & 256 & 0.80 & 3.31E-12 & 6.124 & 2.47E-12 & 6.169 \\
	\cline{3-8}
	 &  & 512 & 0.80 & 2.58E-13 & 3.677 & 2.44E-13 & 3.341 \\
	\cline{1-8} \cline{2-8} \cline{3-8}
	\multirow[t]{5}{*}{6} & \multirow[t]{5}{*}{RK6} & 32 & 0.45 & 4.78E-07 & --- & 4.78E-07 & --- \\
	\cline{3-8}
	 &  & 64 & 0.40 & 3.76E-09 & 6.991 & 3.76E-09 & 6.991 \\
	\cline{3-8}
	 &  & 128 & 0.36 & 2.94E-11 & 6.997 & 2.94E-11 & 6.997 \\
	\cline{3-8}
	 &  & 256 & 0.32 & 8.32E-13 & 5.143 & 8.33E-13 & 5.142 \\
	\cline{3-8}
	 &  & 512 & 0.28 & 4.88E-13 & 0.770 & 4.87E-13 & 0.773 \\
	\cline{1-8} \cline{2-8} \cline{3-8}
	\multirow[t]{5}{*}{7} & \multirow[t]{5}{*}{RK6} & 32 & 0.25 & 6.71E-09 & --- & 5.68E-09 & --- \\
	\cline{3-8}
	 &  & 64 & 0.20 & 2.22E-11 & 8.239 & 1.68E-11 & 8.399 \\
	\cline{3-8}
	 &  & 128 & 0.16 & 8.83E-13 & 4.651 & 8.60E-13 & 4.291 \\
	\cline{3-8}
	 &  & 256 & 0.13 & 1.61E-12 & -0.864 & 1.61E-12 & -0.903 \\
	\cline{3-8}
	 &  & 512 & 0.10 & 4.49E-12 & -1.481 & 4.49E-12 & -1.481 \\
	\cline{1-8} \cline{2-8} \cline{3-8}
	\bottomrule
	\end{tabular}
    \caption{Summary of the results for the 2D advection of a sine wave test using aPrioriMPP and aPrioriT schemes, modified such that the adaptive time-step size is disabled and $C$ is chosen according to (\ref{time_step_modification}). $E^1_\text{GL}$ and $E^1_T$ are the L1 errors of the aPrioriMPP and aPrioriT schemes, respectively. Results are given for varying resolution $N$ and spatial interpolation polynomial degree $p$.}
    \label{tab:convergence}
\end{table}

\subsection{Two-dimensional advection of a square}

We conduct another 2D advection test, keeping $v_x=2$ and $v_y=1$ on the periodic domain $x,y \in [0,1]$, but selecting the discontinuous initial condition:

\begin{equation}
    \overline{\overline{u}}_0 =
    \begin{cases}
    1 & \text{if } 0.25 < x, y < 0.75 \\
    0             & \text{otherwise}.
    \end{cases}
\end{equation} The square is solved up to one period of advection with the aPrioriMPP, aPrioriT, aPosteriori, and aPosterioriB schemes at varying polynomial degree $p$ and with different Runge-Kutta methods.

The maximum principle violations observed from the numerical solutions of these schemes are reported in Table \ref{tab:square2d}. The two-dimensional implementation of aPrioriMPP with the Gauss-Legendre flux quadrature preserves very well the maximum principle with the appropriate SSP Runge-Kutta methods, as well as with RK4. Meanwhile, aPrioriT, which combines \textit{a priori} slope limiting with transverse flux reconstruction, exhibits large violations of the maximum principle regardless of the degree $p$ or the chosen Runge-Kutta method.

aPosteriori and aPosterioriB exhibit maximum principle violations in all cases, with magnitudes that we don't find to depend on $p$, $N$, or $\Delta t$. The violations of the \textit{a posteriori} slope limited schemes are consistently made smaller by using RK4 instead of SSPRK3 and by turning on blending.

\begin{table}[]
	\center
    \begin{tabular}{ll|rr|rr}
    \toprule
     & & \multicolumn{4}{c}{$\delta$} \\
    \cline{3-6}
     $p$ & Integrator & aPrioriMPP & aPrioriT & aPosteriori & aPosterioriB \\
    \midrule
    1 & SSPRK2 & \textbf{-6.79e-11} & -8.82e-04 & -8.72e-03 & -4.14e-03 \\
	\cline{1-6}
	2 & SSPRK3 & \textbf{-9.93e-11} & -1.34e-02 & -1.03e-02 & -1.85e-03 \\
	\cline{1-6}
	\multirow[t]{2}{*}{3} & SSPRK3 & \textbf{-1.72e-12} & -1.34e-02 & -1.17e-02 & -2.09e-03 \\
	& RK4 & \textbf{-5.82e-11} & -1.47e-02 & -1.36e-03 & -3.13e-04 \\
	\cline{1-6}
	\multirow[t]{2}{*}{4} & SSPRK3 & \textbf{-9.43e-12} & -1.61e-02 & -1.30e-02 & -2.05e-03 \\
	& RK4 & \textbf{-8.00e-11} & -1.70e-02 & -2.33e-03 & -2.38e-04 \\
	\cline{1-6}
	\multirow[t]{2}{*}{5} & SSPRK3 & \textbf{-2.71e-11} & -1.61e-02 & -1.66e-02 & -3.43e-03 \\
	& RK4 & \textbf{-7.35e-12} & -1.72e-02 & -3.24e-03 & -4.06e-04 \\
	\cline{1-6}
	\multirow[t]{2}{*}{6} & SSPRK3 & \textbf{-5.69e-11} & -1.72e-02 & -1.48e-02 & -2.97e-03 \\
	& RK4 & \textbf{-2.19e-11} & -1.80e-02 & -6.70e-03 & -1.12e-03 \\
	\cline{1-6}
	\multirow[t]{2}{*}{7} & SSPRK3 & \textbf{-3.42e-11} & -1.72e-02 & -1.65e-02 & -3.41e-03 \\
	& RK4 & \textbf{-1.37e-11} & -1.82e-02 & -7.63e-03 & -1.44e-03 \\
	\cline{1-6}
	\bottomrule
    \end{tabular}
    \caption{Maximum principle violations $\delta$ of aPrioriMPP, aPrioriT, aPosteriori, and aPosterioriB schemes when solving the 2D advection of a square up to $t=1$. The test is repeated for various Runge-Kutta methods and polynomial degrees $p$. Violations smaller in magnitude than -1E-10 are typed in boldface font.}
    \label{tab:square2d}
\end{table}

We also compare the numerical solutions of the various schemes for long time integration using 100 periods. As was the case in 1D, the \textit{a priori} slope-limited schemes use SSPRK3 for $p>2$ since RK4 is observed to result in excess numerical diffusion, while the \textit{a posteriori} slope-limited solutions use RK4 for $p>3$ since it produces smaller violations (seen in Table \ref{tab:square2d} after one period). We observe that the most severe maximum principle violations of our slope-limited schemes occur, if at all, in the first few steps; $\delta$ does not significantly change after one period.

\begin{figure}
    \centering
    \includegraphics[width=\textwidth]{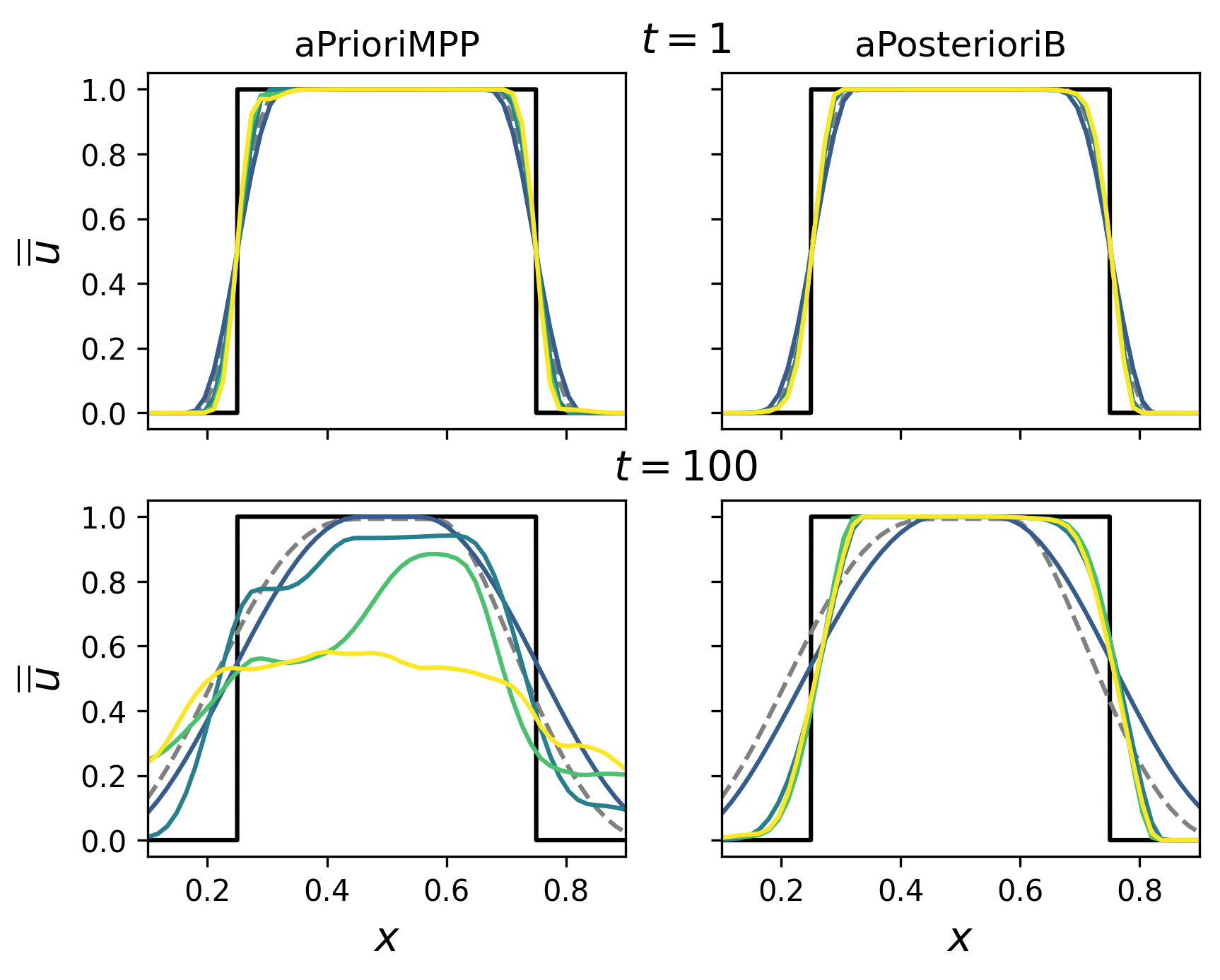}
    \caption{Slices of numerical solutions to the 2D advection of a square along $y=0.5$ at $t=1$ (top row) and $t=100$ (bottom row). Results are shown for aPrioriMPP (left column) and aPosterioriB (right column) schemes and for polynomial degree $p=2$ (dark blue), $p=3$ (light blue), $p=5$ (green), and $p=7$ (yellow). SSPRK2 is used when $p=1$ and SSPRK3 is used when $p=2$. aPosterioriB uses RK4 for $p=3,5,7$ while aPrioriMPP uses SSPRK3 in these cases. The numerical solution of the second-order MUSCL-Hancock is shown in dashed grey for reference. All numerical solutions are presented with a resolution $N=64$.}
    \label{fig:square2d}
\end{figure}

Figure \ref{fig:square2d} shows the numerical solutions of aPrioriMPP and aPosterioriB schemes as well as second-order MUSCL-Hancock after one and one-hundred periods of advection. Neither MUSCL-Hancock nor any of the aPrioriMPP schemes produce violations of the maximum principle. For the shorter time-scale, aPrioriMPP with $p>2$ produces numerical solutions with overall better quality than MUSCL-Hancock despite slight numerical artifacts appearing as $p$ increases. After 100 periods, these numerical artifacts dominate the numerical solution, and the solution quality worsens as $p$ increases for $p=3, 5, 7$. For the discontinuous square, second-order MUSCL-Hancock shows stronger resilience to numerical diffusion than the high-degree, \textit{a priori} slope-limited solutions.

The aPosterioriB schemes, on the other hand, behave much differently. At $p=2$, this \textit{a posteriori} scheme performs similarly to its \textit{a priori} slope-limited counterpart and gives a numerical solution of lower quality than second-order MUSCL-Hancock. Despite that, the $p>2$ aPosterioriB schemes give the highest quality numerical solutions of the schemes presented in Figure \ref{fig:square2d}. These numerical solution hardly change in profile between one and 100 periods, demonstrating an impressive resilience to numerical diffusion.

We clarify that the observed excessive long time-scale numerical diffusion from the high-degree \textit{a priori} limited schemes are a consequence of the chosen resolution $N=64$. By increasing $N$, we can reduce the excessive diffusion of the \textit{a priori} limited schemes. This comes with a significant computational cost, as shown in our later analysis. 

\begin{figure}
    \centering
    \includegraphics[width=\textwidth]{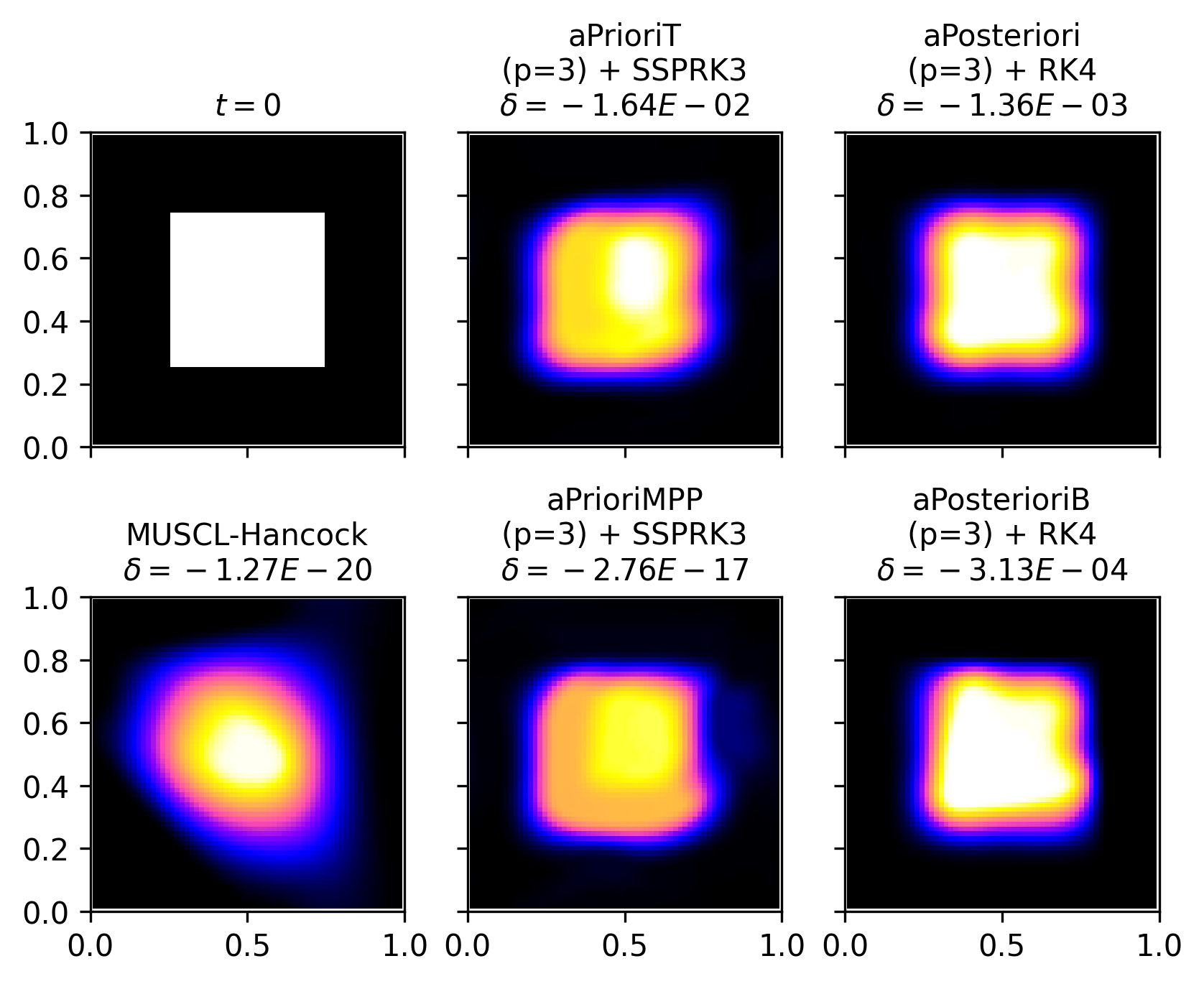}
    \caption{Snapshots of numerical solutions to the 2D advection of a square produced by five different schemes at $t=100$, with the initial condition shown in the top left panel. The value of $\overline{\overline{u}}$ is shaded from 0 to 1 with a color scale varying between black, purple, pink, orange, yellow, and white. Results are shown for a resolution $N=64$ and a polynomial degree $p=3$. The name of the slope limited scheme, chosen Runge-Kutta method, and maximum principle violation $\delta$ are given at the top of each panel.}
    \label{fig:square2_cmap_p3}
\end{figure}

\begin{figure}
    \centering
    \includegraphics[width=\textwidth]{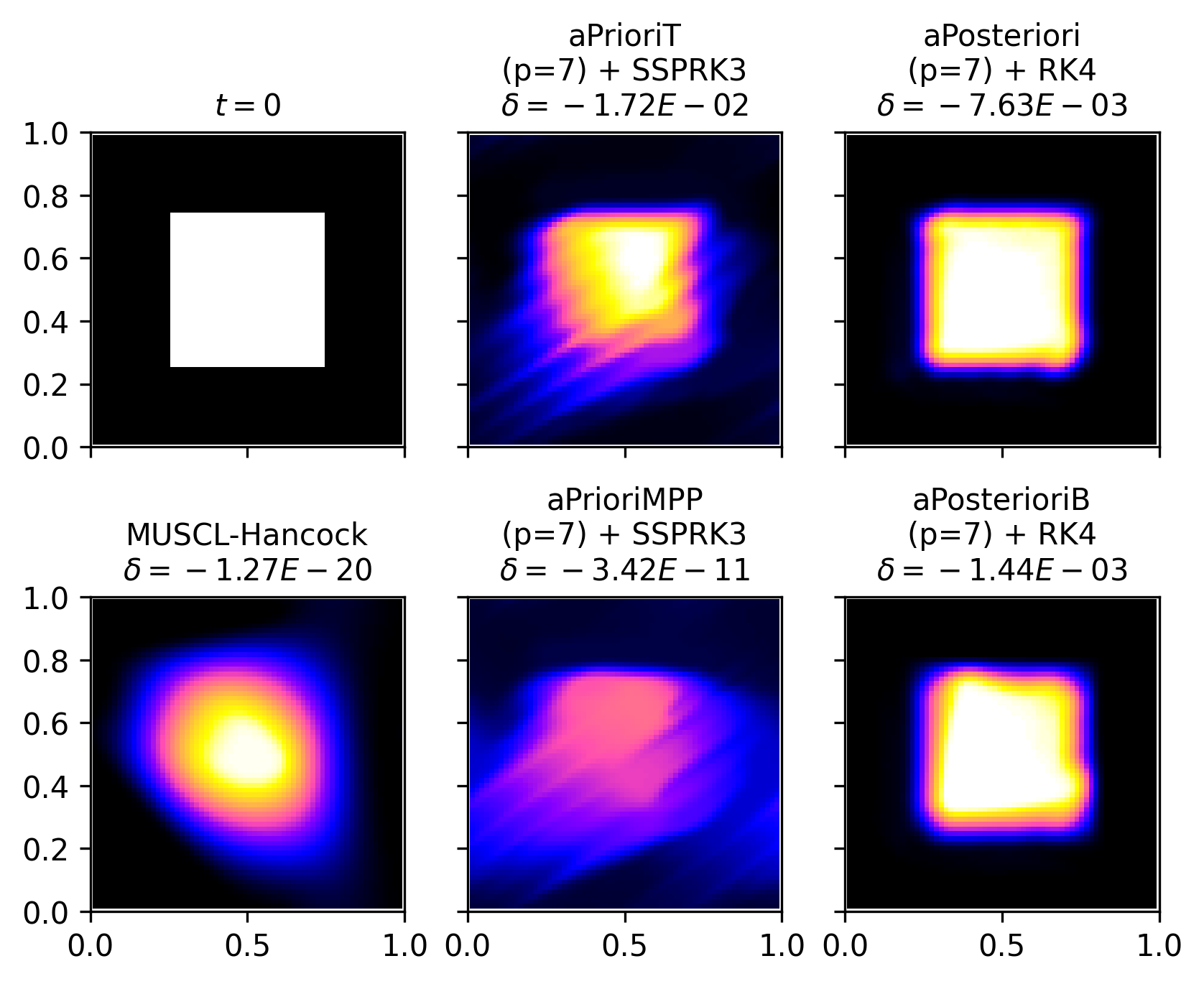}
    \caption{Snapshots of numerical solutions to the 2D advection of a square produced by five different schemes at $t=100$, with the initial condition shown in the top left panel. The value of $\overline{\overline{u}}$ is shaded from 0 to 1 with a color scale varying between black, purple, pink, orange, yellow, and white. Results are shown for a resolution $N=64$ and a polynomial degree $p=7$. The name of the slope limited scheme, chosen Runge-Kutta method, and maximum principle violation $\delta$ are given at the top of each panel.}
    \label{fig:square2_cmap_p7}
\end{figure}

Figures \ref{fig:square2_cmap_p3} and \ref{fig:square2_cmap_p7} show color maps of the numerical solutions of all four high-order schemes schemes at $t=100$. with $p=3$ and $p=7$, respectively. MUSCL-Hancock is also shown in both figures. When comparing the two MPP schemes in Figure \ref{fig:square2_cmap_p3}, MUSCL-Hancock and aPrioriMPP ($p=3$), it is evident that aPrioriMPP better preserves the sharp gradients of the square at its edges, despite the visibility of some numerical artifacts. When $p=7$, we again see that the aPrioriMPP numerical solution is dominated by numerical diffusion and artifacts and is lower in quality than the MUSCL-Hancock solution. This degradation is reduced somewhat by using transverse flux reconstruction instead of the Gauss-Legendre quadrature, as is the case for aPrioriT, but this variation of scheme still exhibits large violations of the maximum principle.

Both aPosteriori and aPosterioriB show great solution quality compared to MUSCL-Hancock and the \textit{a priori} slope limited schemes, giving generally better results as $p$ increases from 3 to 7. Aside from their maximum principle violations, both \textit{a posteriori} slope limiting schemes greatly outperform MUSCL-Hancock in terms of solution quality. In both cases, aPosterioriB exhibits smaller maximum principle violations while also showing more numerical artifacts along the leading edge of the square.

\subsection{Rotation of a slotted cylinder}

We include the classical test of the two-dimensional rotation of a slotted disk with the initial condition

\begin{equation}
\overline{\overline{u}}_0(x,y) =
    \begin{cases}
        1 & x^2 + (y-0.5)^2 < 0.3^2, \quad |x| > 0.025, \quad y > 0.7\\
        0 & \text{otherwise}
    \end{cases}
\end{equation} defined on $x,y \in [-1,1]$ with a Dirichlet boundary fixed at 0 and a non-uniform velocity field $\langle v_x, v_y \rangle = \langle -y, x \rangle$. This test introduces additional complexity since the velocity is not uniform. 

We again report the maximum principle violations of aPrioriMPP, aPrioriT, aPosteriori, and aPosterioriB schemes after one period of rotation for different polynomial degrees $p$ and Runge-Kutta methods. These results are provided in Table \ref{tab:slottedcylinder}. aPrioriMPP gives very good MPP results at every $p$ from 1 to 7 and for the SSPRK2, SSPRK3, and RK4 Runge-Kutta methods. aPrioriT exhibits large violations in all cases. The \textit{a posteriori} slope limited schemes violate the maximum principle in every case, but with a magnitude that is consistently reduced by using RK4 and blending the revised fluxes. For this experiment, the magnitude of the \textit{a posteriori} slope-limiting violations do not strictly decrease by increasing $p$ or $N$, nor by decreasing the time-step size.

\begin{table}[]
	\center
    \begin{tabular}{ll|rr|rr}
    \toprule
     & & \multicolumn{4}{c}{$\delta$} \\
    \cline{3-6}
     $p$ & Integrator & aPrioriMPP & aPrioriT & aPosteriori & aPosterioriB \\
    \midrule
    1 & SSPRK2 & \textbf{-6.79e-11} & -8.82e-04 & -8.72e-03 & -4.14e-03 \\
	\cline{1-6}
	2 & SSPRK3 & \textbf{-9.93e-11} & -1.34e-02 & -1.03e-02 & -1.85e-03 \\
	\cline{1-6}
	\multirow[t]{2}{*}{3} & SSPRK3 & \textbf{-1.72e-12} & -1.34e-02 & -1.17e-02 & -2.09e-03 \\
	 & RK4 & \textbf{-5.82e-11} & -1.47e-02 & -1.36e-03 & -3.13e-04 \\
	\cline{1-6}
	\multirow[t]{2}{*}{4} & SSPRK3 & \textbf{-9.43e-12} & -1.61e-02 & -1.30e-02 & -2.05e-03 \\
	 & RK4 & \textbf{-8.00e-11} & -1.70e-02 & -2.33e-03 & -2.38e-04 \\
	\cline{1-6}
	\multirow[t]{2}{*}{5} & SSPRK3 & \textbf{-2.71e-11} & -1.61e-02 & -1.66e-02 & -3.43e-03 \\
	 & RK4 & \textbf{-7.35e-12} & -1.72e-02 & -3.24e-03 & -4.06e-04 \\
	\cline{1-6}
	\multirow[t]{2}{*}{6} & SSPRK3 & \textbf{-5.69e-11} & -1.72e-02 & -1.48e-02 & -2.97e-03 \\
	 & RK4 & \textbf{-2.19e-11} & -1.80e-02 & -6.70e-03 & -1.12e-03 \\
	\cline{1-6}
	\multirow[t]{2}{*}{7} & SSPRK3 & \textbf{-3.42e-11} & -1.72e-02 & -1.65e-02 & -3.41e-03 \\
	 & RK4 & \textbf{-1.37e-11} & -1.82e-02 & -7.63e-03 & -1.44e-03 \\
	\cline{1-6}	\bottomrule
    \end{tabular}
    \caption{Maximum principle violations $\delta$ of aPrioriMPP, aPrioriT, aPosteriori, and aPosterioriB schemes when solving the rotation of a slotted disk up to $t=2\pi$. The test is repeated for various Runge-Kutta methods and polynomial degrees $p$. Violations smaller in magnitude than -1E-10 are typed in boldface font.}
    \label{tab:slottedcylinder}
\end{table}

The slotted disk is also solved up to 10 periods of advection to explore the effect of long time integration. Our usual policy is implemented wherein \textit{a priori} slope-limited schemes do not use higher-than-third order time integration (SSPRK3) while \textit{a posteriori} slope-limited schemes use up to fourth-order time integration (RK4). 

\begin{figure}
    \centering
    \includegraphics[width=\textwidth]{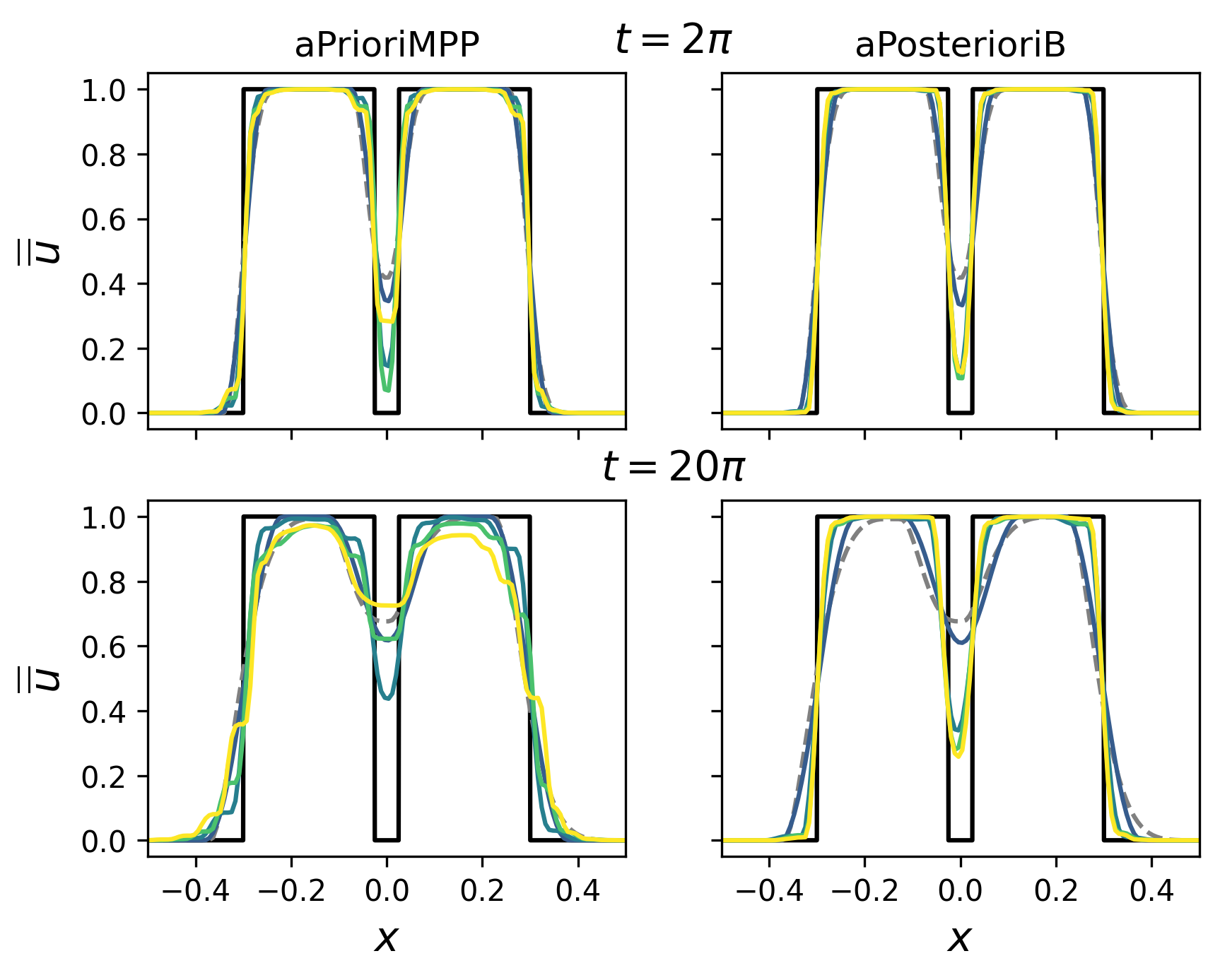}
    \caption{Slices of numerical solutions to the rotation of a slotted cylinder along $y=0.5$ at $t=2\pi$ (top row) and $t=20\pi$ (bottom row). Results are shown for aPrioriMPP (left column) and aPosterioriB (right column) schemes and for polynomial degree $p=2$ (dark blue), $p=3$ (light blue), $p=5$ (green), and $p=7$ (yellow). SSPRK2 is used when $p=1$ and SSPRK3 is used when $p=2$. aPosterioriB uses RK4 for $p=3,5,7$ while aPrioriMPP uses SSPRK3 in these cases. The numerical solution of second-order MUSCL-Hancock is shown in dashed grey for reference. All numerical solutions are presented with a resolution $N=256$.}
    \label{fig:slottedcylinder}
\end{figure}

The numerical solutions of aPrioriMPP and aPosterioriB schemes with increasing polynomial degree $p$, as well as MUSCL-Hancock, are shown after one and ten periods of advection in Figure \ref{fig:slottedcylinder}. After one period, the \textit{a priori} and \textit{a posteriori} slope limited schemes at the same $p$ give numerical solutions that look quite similar, with slightly more numerical artifacts in the higher-$p$ aPrioriMPP solutions. At ten periods, the quality of aPrioriMPP significantly worsens as $p$ increases after around $p=3$. aPosterioriB, on the other hand, gives generally better quality solutions as $p$ increases. Both high-order slope limited schemes at both time-scales outperform second-order MUSCL-Hancock in terms of numerical diffusion. 

\begin{figure}
    \centering
    \includegraphics[width=\textwidth]{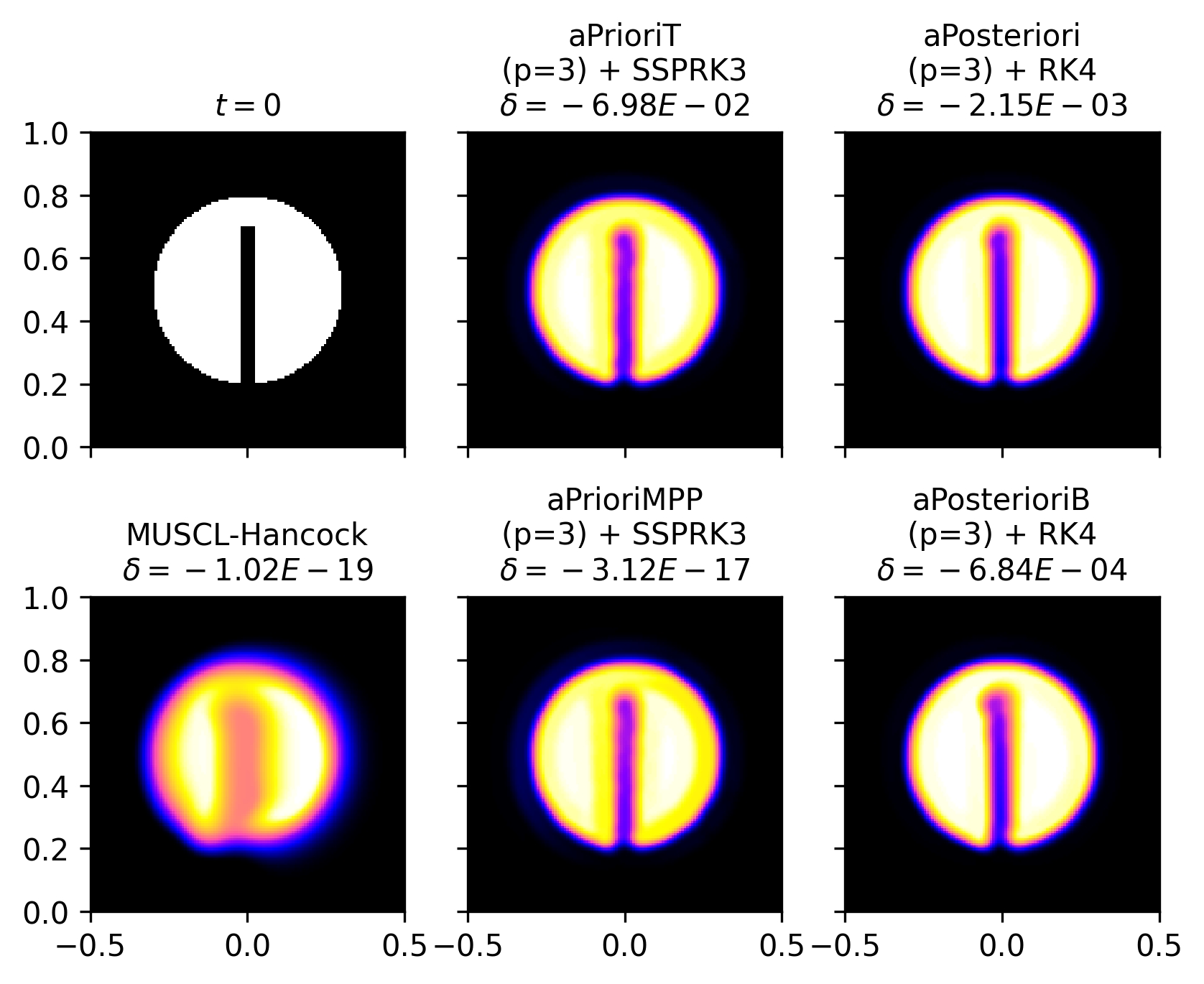}
    \caption{Snapshots of numerical solutions to the rotation of a slotted cylinder produced by five different schemes at $t=20\pi$, with the initial condition shown in the top left panel. The value of $\overline{\overline{u}}$ is shaded from 0 to 1 with a color scale varying between black, purple, pink, orange, yellow, and white. Results are shown for a resolution $N=256$ and a polynomial degree $p=3$. The name of the slope limited scheme, chosen Runge-Kutta method, and maximum principle violation $\delta$ are given at the top of each panel.}
    \label{fig:slottedcylinder_cmap_p3}
\end{figure}

\begin{figure}
    \centering
    \includegraphics[width=\textwidth]{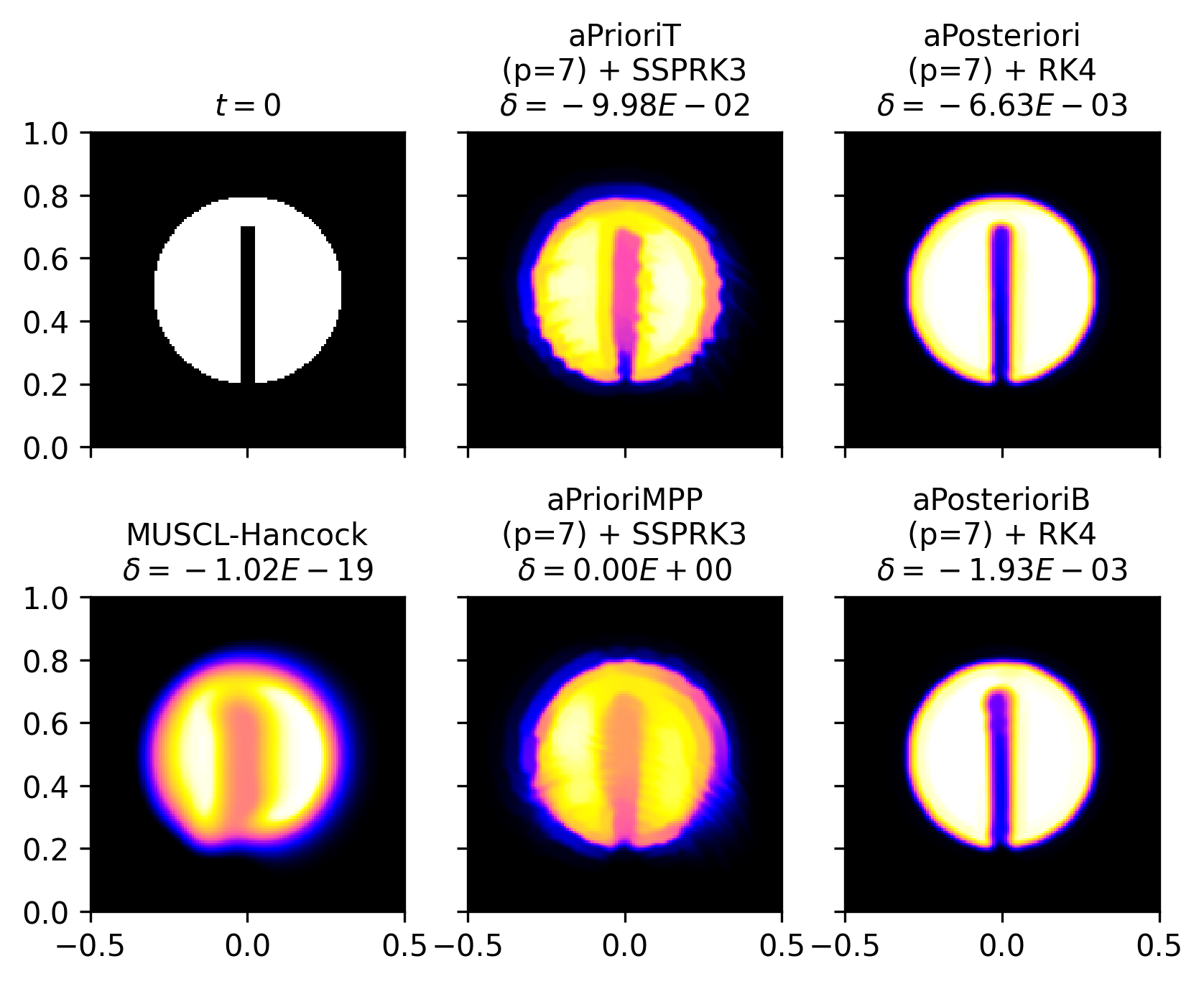}
    \caption{Snapshots of numerical solutions to the rotation of a slotted cylinder produced by five different schemes at $t=20\pi$, with the initial condition shown in the top left panel. The value of $\overline{\overline{u}}$ is shaded from 0 to 1 with a color scale varying between black, purple, pink, orange, yellow, and white. Results are shown for a resolution $N=256$ and a polynomial degree $p=7$. The name of the slope limited scheme, chosen Runge-Kutta method, and maximum principle violation $\delta$ are given at the top of each panel.}
    \label{fig:slottedcylinder_cmap_p7}
\end{figure}

The numerical solutions of aPrioriMPP, aPrioriT, aPosteriori, and aPosterioriB are compared at $t=20\pi$ with $p=3$ in Figure \ref{fig:slottedcylinder_cmap_p3} and $p=7$ in Figure \ref{fig:slottedcylinder_cmap_p7}. At $p=3$, aPrioriMPP and aPrioriT appear similar, with aPrioriMPP approximately preserving the maximum principle and aPrioriT causing large violations. The \textit{a posteriori} slope-limited schemes give slightly less diffused results, with aPosterioriB giving a smaller violation of the maximum principle as well as more apparent numerical artifacts. This is also seen for aPosteriori and aPosterioriB at $p=7$. At $p=7$, the numerical artifacts in the aPrioriMPP and aPrioriT solutions are strong, while aPosteriori and aPosterioriB give improved results from the $p=3$ case. All four high-order schemes at both $p=3$ and $p=7$ exhibit less numerical diffusion than MUSCL-Hancock at the same resolution. As was the case with the 2D square experiment, the selection of $N$ ultimately determines the amount of numerical diffusion produced by the high-order schemes.

\subsection{Cost analysis}

The computational cost of a finite volume scheme depends largely on the choice of flux integral method. The transverse flux reconstruction requires the nodal reconstruction of only one point per cell face, while the Gauss-Legendre quadrature requires $\ceil{(p+1)/2}$ points. Each reconstructed node depends on a large matrix multiplication, and each node corresponds to a Riemann problem, which can be the most expensive step in the solution for non-linear conservation laws. The transverse reconstruction presents a clear advantage in terms of cost. However, it was shown to result in large maximum principle violations when used in conjunction with the \textit{a priori} limiting method, so, in our opinion, only the \textit{a posteriori} slope limiting method is allowed to benefit from the low-cost flux reconstruction.

GPUs, which are optimized for large matrix multiplication problems, could potentially mitigate the cost of the reconstructions required by the Gauss-Legendre quadrature. Therefore, we run this experiment on both CPU and GPU implementations. The CPU version of the code is written in standard Python and uses the NumPy library for optimized array multiplications. The GPU version includes the CuPy library which allows for these same multiplications to occur on a GPU. An NVIDIA A100 is used for this experiment and no parallelization over multiple GPUs was implemented.

A timing comparison was conducted for the 2D advection of a square for two different slope limiting methods: aPrioriMPP, the \textit{a priori} slope limiting method using the Gauss-Legendre quadrature, and aPosterioriB, the \textit{a posteriori} slope limiting method using transverse flux reconstruction. The adaptive time-step size of aPrioriMPP is turned off for this experiment and RK6 is used for the schemes with $p>3$. Computational speed is reported as the number of cells in each Runge-Kutta stage updated per second, taken as an average over ten time-steps. Since the compute time increases linearly with the number of Runge-Kutta stages, variations in speed between Runge-Kutta methods are not visible with this metric; only the variation in speed of the different spatial discretization methods due to changes in $p$ are visible.

\begin{figure}
    \centering
    \includegraphics[width=\textwidth]{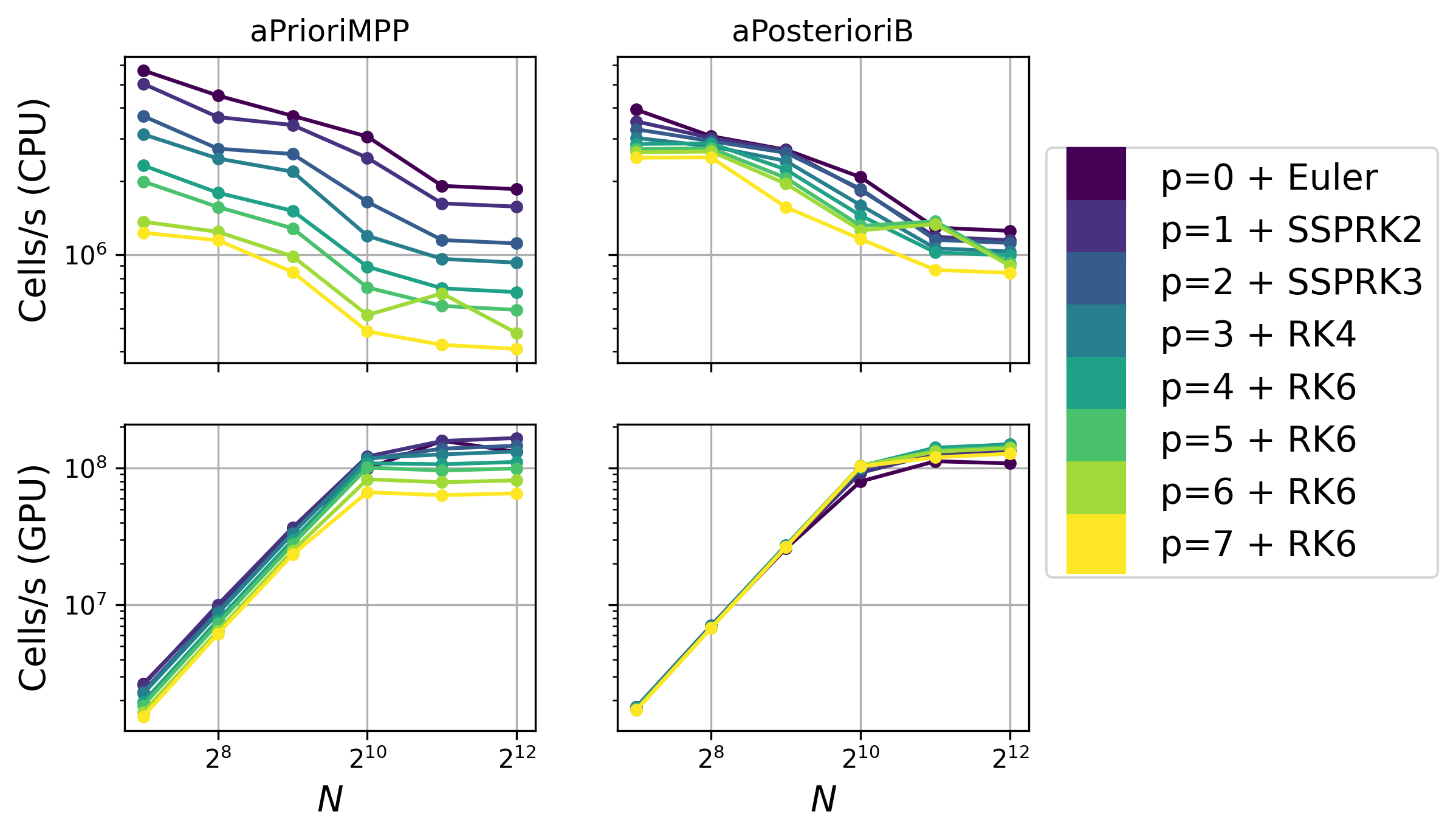}
    \caption{The compute speed measured for aPrioriMPP (left column) and aPosteriori (right column) schemes and for varying polynomial degree $p$ and Runge-Kutta method. The adaptive time-step size of aPrioriMPP is turned off for this experiment. Compute speed is measured as the number of cells updated in each Runge-Kutta stage per second, taken as an average over ten time-steps. Results are shown for implementations on both CPU (top row) and GPU (bottom row) platforms.}
    \label{fig:timing}
\end{figure}

The results of this experiment are depicted in Figure \ref{fig:timing}. On the CPU, both the aPriori and aPosteriori schemes show an increase in cost with an increasing polynomial degree $p$. However, the cost of the \textit{a priori} limiting method escalates at a much faster rate with $p$ than that of the \textit{a posteriori} limiting method due to the greater number of nodal reconstructions required by the Gauss-Legendre quadrature. It's worth noting that the timing of even/odd pairs of the aPrioriMPP schemes—$p=0,1$, $p=2,3$, $p=4,5$, and $p=6,7$—appears to cluster together because the members of each pair use the same number of quadrature points. This clustering effect is not observed for the aPosterioriB schemes, which utilize a transverse flux reconstruction. Furthermore, the per-cell cost of both methods on the CPU increases as the size $N$ grows larger, as the computational overhead of handling large arrays becomes more significant.

Conversely, on the GPU, the per-cell cost of both methods \textit{decreases} with the problem size until the threads of the GPU become saturated at about $N=2^{11}$. At this size, the schemes are about two orders of magnitudes faster on GPU than on CPU. The disparity in cost between the \textit{a priori} and \textit{a posteriori} schemes is notably reduced on the GPU, highlighting the GPU's efficiency in intensive matrix multiplication processes. Likewise, aPosterioriB does not become significantly more expensive with increasing $p$ on the GPU due to the reduced number of required matrix multiplications from the transverse flux reconstruction.

Even though the \textit{a posteriori} limiting method avoids the expensive Gauss-Legendre quadrature at each face required by the \textit{a priori} limiting method, it still comes with a significant cost of its own: the fallback scheme, which includes detecting troubled cells and revising their fluxes. At low spatial degree ($p=3$), we find that the fallback scheme accounts for roughly 1/2 of the computational time while at a higher degree ($p=7$), it accounts for 1/3.

\section{Discussion}
\label{sec:discussion}

The results obtained from our numerical tests provide strong validation for the accuracy of our high-order finite volume methods and smooth extrema detection implementation. The advection test featuring a smooth sine wave showcases the high-order capabilities of our base finite volume scheme in two spatial dimensions. Notably, these results were achieved with \textit{a priori} slope limiting enabled, confirming that our smooth extrema detection correctly disables slope limiting in smooth regions. Moreover, this success was achieved using either the Gauss-Legendre or transverse flux reconstructions. Nonetheless, lower errors were found with the transverse flux reconstruction at some even orders of accuracy.

In all conducted experiments, our implementation of Zhang \& Shu's \textit{a priori} slope limiting method \cite{zhang2011maximum, zhang2010maximum} consistently demonstrated the absence of maximum principle violations when paired with forward Euler, SSPRK2, or SSPRK3. Recall that the latter two methods are equivalent to a convex combination of forward Euler steps, making them Strong Stability Preserving (SSP), and that Zhang \& Shu's spatial discretization is guaranteed to preserve a strict maximum principle when solved with such Runge-Kutta methods \cite{zhang2011maximum, zhang2010maximum}.

RK4, on the other hand, is not SSP. However, it consistently demonstrated quasi-SSP behavior; For every scheme in our experiments that was MPP when solved with an SSP method, it was also MPP when solved with RK4. This quasi-SSP behavior has been described for RK4 (see \cite{sanz2010positivity}), but not for other non-SSP Runge-Kutta methods such as RK6. In fact we did not find RK6 to enable our slope-limited schemes to avoid maximum principle violations. Despite being quasi-SSP, RK4 is found to result in slightly more numerical diffusion than SSPRK3 for the \textit{a priori} slope-limited schemes with large values of $p$. We observe the opposite effect for the \textit{a posteriori} slope-limited schemes.

We implemented modifications to the Zhang \& Shu \textit{a priori} slope limited scheme \cite{zhang2011maximum, zhang2010maximum} that maintain the Maximum Principle Preserving (MPP) property while reducing computational cost and complexity. For instance, the centroid set of points, employed in both one- and two-dimensional cases, proves to be a suitable substitute for the Gauss-Lobatto quadrature points when computing the \textit{a priori} slope limiter $\theta_{ij}$. Additionally, our experiments revealed that the reduced CFL factor $C_\text{MPP}$ is unnecessary when using an adaptive time-step size of Huang \textit{et al} \cite{huang2023high}.

In the 1D experiments, we determined that the \textit{a priori} and \textit{a posteriori} limited schemes give similar, high-quality numerical solutions when implemented at high-order. This is true even for long time integration. This is not the case for the 2D experiments, where we see a difference between the results after one period of advection and after many periods.

After a single period of advection, both the \textit{a priori} and \textit{a posteriori} limited solutions to the 2D problems exhibit good solution quality for all of the third-or-higher-order schemes shown ($p\ge2$) across all experiments featuring discontinuous initial conditions. Notably, the numerical solutions of second-order MUSCL-Hancock performed quite well for small time integration, outperforming our $p=2$ schemes for the composite profile and the discontinuous 2D square.

After many periods of advection (100 for the composite profile and discontinuous square and 10 for the slotted cylinder), we see a trade-off between the quality of high-order solutions and the magnitude of their maximum principle violations. For instance, the results of the aPrioriMPP schemes produce no maximum principle violations, but their numerical solutions suffer from significant diffusion as $p$ increases after about $p=4$. At the largest value $p=7$, the numerical solution of aPrioriMPP becomes dominated by oscillatory artifacts. The worst example of this is seen in the 2D advection of the square.

It was clarified that the excessive diffusion produced by the high-degree aPrioriMPP schemes in the 2D experiments are the consequence of the resolution $N$ chosen to demonstrate its contrast with the corresponding \textit{a posteriori} limiting schemes. This numerical diffusion can always be reduced by increasing $N$, but this comes with a particularly large computational cost for the high-order, \textit{a priori} limited schemes.

So, to use an aPrioriMPP scheme, for problems with long time integration, we are compelled to choose the sweet spot value of $p=3$ or $p=4$. On the other hand, the \textit{a posteriori} limited schemes typically improve in quality at long time-scales as we increase $p$. Of course, this comes at the cost of higher maximum principle violations.

An additional caveat with the \textit{a priori} limiting method is that it produces large violations (>1\%) when used with transverse flux reconstruction. This observation is not unique to our work; McCorquodale \& Colella \cite{mccorquodale2011high} use an \textit{a priori} limiter as well as transverse flux reconstruction and similar maximum principle violations on the order of >1\% are observed when they advect a discontinuous square. Their limiter is not the Zhang \& Shu limiter \cite{zhang2011maximum, zhang2010maximum}, but is still conceptually similar.

We demonstrated that the magnitude of violations from \textit{a posteriori} limited schemes can be reduced by applying Vilar and Abgrall's convex blending of the corrected fluxes \cite{vilar2024posteriori}. Additionally, the magnitude of these violations were consistently reduced when RK4 is used instead of SSPRK3. Interestingly, the \textit{a posteriori} limited schemes that produce smaller violations of the maximum principle tend to exhibit greater numerical diffusion and more pronounced numerical artifacts. This highlights a fundamental trade-off where reducing maximum principle violations often comes at the cost of introducing other undesirable characteristics in the numerical solution.

In the case of strict positivity preservation, the magnitude of the maximum principle violation of a numerical scheme limits the dynamical range that can be captured for a given scalar. For example, suppose we are solving Euler's equation.  Our solver will fail for any value of the mass density $\rho<0$. If we use an \textit{a posteriori} slope limited scheme for which we expect the maximum principle violation to satisfy $|\delta| \in (0,0.001)$, the scheme will certainly fail if the density contrast is larger than $10^{3}$ since maximum principle violations will produce negative densities. However, if the density contrast remains smaller than $10^{3}$, the magnitude of the maximum principle violations of our scheme will result in the positivity of the solution $\rho \ge 0$. Thus, a smaller $|\delta|$ of a given scheme allows for higher density contrasts.

While it is feasible to implement the \textit{a priori} limiting method at high order and maintain a strict maximum principle, we demonstrate that this requires a more computationally expensive flux quadrature. Specifically, we find that the transverse reconstruction incurred significantly lower computational costs on a CPU. However, it is incompatible with \textit{a priori} schemes due to substantial maximum principle violations resulting from this combination. On a GPU, the difference in cost between the two quadratures was mostly mitigated. This is because the additional matrix multiplications associated with the Gauss-Legendre quadrature, which contribute to its higher cost on a CPU, become relatively insignificant on a GPU due to its optimized performance for such tasks. We point out that a GPU implementation must solve problems of a sufficiently large size for the increase in speed to be worthwhile. For instance, our GPU code was only faster than the purely CPU code for roughly $N > 256$ in 2D. This is due to the overhead cost of communicating arrays between the CPU and GPU, which is only insignificant when the arrays are large. 

\section{Conclusion}
\label{sec:conclusion}

In this study, we develop a novel and experimental approach to compare \textit{a priori} and \textit{a posteriori} slope limiters for high-order finite volume schemes. Our \textit{a priori} slope limiting schemes are based on the method developed by Zhang \& Shu \cite{zhang2011maximum,zhang2010maximum}, while our \textit{a posteriori} schemes follow a flux revision procedure with a MUSCL fallback scheme \cite{velasco2023spectral}. To assess the relative performance of the two types of schemes, the linear-advection equation is solved in one- and two-dimensions for various benchmark problems at various spatial polynomial degree $p$. In these experiments, the  schemes are compared based on the following figures of merit:
	
	\begin{itemize}
	  \item Ability to preserve the maximum principle or the positivity of the solution.
	  \item Numerical solution quality for long time integration.
	  \item Computational cost as it scales with the resolution $N$ and spatial polynomial degree $p$.
	\end{itemize}

In the one-dimensional case, we find that the implementations of the \textit{a priori} limited schemes for $p>2$ achieve all three goals quite well. The \textit{a posteriori} limited schemes result in similar quality, but consistently produce (small) maximum principle violations. We observe that RK4, despite not being a Strong Stability Preserving method, also strictly preserves the maximum principle when used with the \textit{a priori} limited schemes in one- and two-dimensions.

In two dimensions, the conclusions are more nuanced due to the significant impact of computational cost, among other factors. The \textit{a priori} limiting schemes lead to large maximum principle violations when using the cost-effective transverse flux reconstruction. This necessitates the use of the more expensive Gauss-Legendre quadrature, with a number of nodal reconstructions that grows linearly with $p$. On the other hand, the numerical solution of the \textit{a posteriori} limiting schemes do not exhibit a significant change between these flux reconstruction methods. Thus, these schemes can benefit from the low-cost transverse flux reconstruction, making them highly competitive against their \textit{a priori} counterparts in terms of speed. This speed comparison is more relevant to implementations on classic CPU architecture; Timing experiments revealed that, for high enough resolutions, the cost difference between the two schemes is dramatically reduced when the computations are performed on GPUs. 

The \textit{a priori} schemes also suffer from another drawback: their solution quality deteriorates for long time integration as $p$ increases beyond approximately $p=4$. When $p>4$ and there isn't enough cell resolution, numerical artifacts from the high-degree interpolation polynomials begin to dominate the solution, resulting in excessive diffusion. In contrast, the \textit{a posteriori} limited schemes exhibit remarkable resilience to numerical diffusion and generally yield better results as $p$ increases, taking full advantage of the steep gradients offered by the high-degree polynomials. We observe that at the same, modest resolution, the \textit{a posteriori} scheme excels at long time-scales while the \textit{a priori} scheme suffers from diffusion.

While the \textit{a posteriori} schemes outperform their \textit{a priori} counterparts in terms of solution quality for long time integration and cost, they still have the issue of maximum principle violations. In fact, the violations of the \textit{a posteriori} schemes are typically greater in the two-dimensional than in the one-dimensional case. These violations persist irrespective of time-size $\Delta t$, $N$, or $p$. While the convex blending of revised fluxes proposed by Vilar \& Abgrall \cite{vilar2024posteriori} can reduce these violations, it comes at the expense of a slight decrease in solution quality, mirroring the degradation observed in \textit{a priori} limited solutions.

Our study indicates that \textit{a posteriori} limiting is an excellent choice for problems where some amount of maximum principle violations can be tolerated. The combination of long time-scale solution quality and cost efficiency offered by these schemes are not matched by second-order MUSCL-Hancock or the \textit{a priori} slope-limited schemes. In the case where maximum principle violations must be strictly bounded but long time-scale diffusion must be avoided, \textit{a priori} slope-limited schemes are a suitable, albeit more expensive option. These schemes also come with the caveat that if they are implemented at too high of order ($p>4$) or with too low of a resolution $N$, their numerical solutions might become dominated by artifacts for long integration times.

Future research directions could focus on exploring further modifications to \textit{a posteriori} slope limiting methods aimed at reducing or preventing maximum principle violations entirely. Additionally, extending this work beyond the linear advection equation to encompass non-linear conservation laws, such as Euler's equation, would provide valuable insights into the performance and applicability of these slope limiting schemes in broader contexts.

%\appendix
%\section{Example Appendix Section}
%\label{app1}

%% BIBLIOGRAPHY

\bibliographystyle{elsarticle-num} 
\bibliography{manuscript}



%\begin{thebibliography}{00}
%
%%% For numbered reference style
%%% Text of bibliographic item
%
%
%\bibitem{lamport94}
%  \textit{\LaTeX: a document preparation system},
%  Addison Wesley, Massachusetts,
%\end{thebibliography}

\end{document}